\documentclass{amsart}

\usepackage[all]{xy}
\usepackage{amssymb}

\theoremstyle{plain}
\newtheorem{theorem}{Theorem}[section]
\newtheorem{proposition}[theorem]{Proposition}
\newtheorem{lemma}[theorem]{Lemma}
\newtheorem{corollary}[theorem]{Corollary}

\theoremstyle{definition}

\newtheorem{definition}[theorem]{Definition}

\newtheorem{remark}[theorem]{Remark}

\renewcommand{\AA}{\mathbb{A}}
\newcommand{\BB}{\mathbb{B}}

\newcommand{\ZZ}{\mathbb{Z}}

\newcommand{\calD}{\mathcal{D}}
\newcommand{\calE}{\mathcal{E}}

\newcommand{\calO}{\mathcal{O}}
\newcommand{\calP}{\mathcal{P}}

\newcommand{\calS}{\mathcal{S}}
\newcommand{\calT}{\mathcal{T}}
\newcommand{\calV}{\mathcal{V}}

\newcommand{\Aut}{\mathrm{Aut}}
\newcommand{\Gr}{\mathrm{Gr}}
\newcommand{\Ho}{\mathrm{H}}
\newcommand{\Hom}{\mathrm{Hom}}
\newcommand{\id}{\mathrm{id}}
\newcommand{\iso}{\mathrm{iso}}

\newcommand{\objects}{\mathrm{obj}}

\newcommand{\rmb}{\mathrm{b}}

\begin{document}

\title{Determinant functors on triangulated categories}
\author{Manuel Breuning}
\address{Department of Mathematics, King's College London, Strand,
London WC2R 2LS, United Kingdom} \email{manuel.breuning@kcl.ac.uk}
\thanks{The author was supported by EPSRC grant GR/S91772/01.}
\date{October 13, 2006}

\begin{abstract}
We study determinant functors which are defined on a triangulated
category and take values in a Picard category. The two main results
are the existence of a universal determinant functor for every small
triangulated category, and a comparison theorem for determinant
functors on a triangulated category with a non-degenerate bounded
t-structure and determinant functors on its heart. For a small
triangulated category $\calT$ we give a natural definition of groups
$K_0(\calT)$ and $K_1(\calT)$ in terms of the universal determinant
functor on $\calT$, and we show that $K_i(\calT)\cong K_i(\calE)$
for $i=0$ and $1$ if $\calT$ has a non-degenerate bounded
t-structure with heart $\calE$.
\end{abstract}

\maketitle


\section{Introduction}

The systematic study of determinant functors began with the
construction of a determinant functor on the category of perfect
complexes of $\calO_X$-modules on a scheme $X$ by Knudsen and
Mumford \cite{KnusemMumford76}. For every perfect complex $A$ they
defined a graded invertible $\calO_X$-module $f(A)$ (which is
essentially the alternating product of the top exterior powers of a
bounded locally free resolution of $A$), and they showed that $f$
can be made into a functor if one only considers quasi-isomorphisms
of perfect complexes. Furthermore for every short exact sequence
$0\to A\to B\to C\to 0$ of perfect complexes they defined an
isomorphism $f(B)\cong f(A)\otimes f(C)$ in the category of graded
invertible $\calO_X$-modules which satisfies a list of natural
properties.

The functor $f$ factors through the category whose objects are
perfect complexes of $\calO_X$-modules and whose morphisms are
isomorphisms in the derived category of $\calO_X$-modules. However
it was observed by Knudsen and Mumford \cite[p.\
44]{KnusemMumford76} that in general it is not possible to define a
natural isomorphism $f(B)\cong f(A)\otimes f(C)$ for every
distinguished triangle $A\to B\to C\to TA$ of perfect complexes.

The results of Knudsen and Mumford were extended to a theory of
determinant functors on exact categories with values in
(commutative) Picard categories (see Deligne \cite{Deligne87},
Knudsen \cite{Knudsen02}). These determinant functors have found
applications in various areas of mathematics. However since one is
often interested in applying determinant functors to objects in a
derived category, one would like to have determinant functors which
are not only well-behaved on short exact sequences but also take
into account the triangulated structure of these categories. The
purpose of this paper is to suggest a definition of determinant
functors on triangulated categories and to demonstrate that this
definition leads to an interesting theory.

Recall that a Picard category $\calP$ is a category with a bifunctor
$\otimes$ and associativity and commutativity isomorphisms
satisfying several natural properties. We define a determinant
functor $f$ on a triangulated category $\calT$ with values in a
Picard category $\calP$, in symbols $f:\calT\to\calP$, to be a pair
$f=(f_1,f_2)$ consisting of a functor $f_1: \calT_{\iso}\to\calP$
(where $\calT_{\iso}$ is the category obtained from $\calT$ by
restricting morphisms to isomorphisms) together with an isomorphism
$f_2(\Delta): f_1(B)\cong f_1(A)\otimes f_1(C)$ for every
distinguished triangle $\Delta: A\to B\to C\to TA$ in $\calT$, such
that $f_2$ is natural in isomorphisms of distinguished triangles and
compatible with direct sums and octahedral diagrams in $\calT$. We
will also define morphisms of determinant functors and therefore
obtain a category $\det(\calT,\calP)$ of determinant functors on
$\calT$ with values in $\calP$.

A determinant functor $f:\calT\to\calV$ is called universal if for
every Picard category $\calP$ the induced functor
\begin{equation*}
\Hom^\otimes(\calV,\calP)\to\det(\calT,\calP)
\end{equation*}
which sends $M$ to $M\circ f$ is an equivalence of categories. Here
$\Hom^\otimes(\calV,\calP)$ denotes the category of monoidal
functors from $\calV$ to $\calP$, and $M\circ f$ is the natural
composition of a determinant functor and a monoidal functor. Every
small triangulated category $\calT$ has a universal determinant
functor $f:\calT\to\calV$. We will prove this by giving an explicit
construction of $\calV$ in terms of generators and relations.

The existence of a universal determinant functor $f: \calE\to\calV$
on a small exact category $\calE$ was shown by Deligne \cite[\S
4]{Deligne87}. His proof shows that for $i=0$ and $1$ the Quillen
$K$-group $K_i(\calE)$ is canonically isomorphic to the homotopy
group $\pi_i(\calV)$. For a small triangulated category $\calT$ this
motivates the definition $K_i(\calT)=\pi_i(\calV)$ for $i=0, 1$
where $f: \calT\to\calV$ is any universal determinant functor on
$\calT$. It follows easily from the definition of a universal
determinant functor that these $K$-groups are well-defined and
functorial. Furthermore it is not difficult to see that $K_0(\calT)$
is the Grothendieck group of $\calT$.

Next we assume that $\calT$ is a triangulated category with a
non-degenerate bounded t-structure and we let $\calE$ be its heart.
We will show that for every Picard category $\calP$ there is a
natural equivalence of categories
\begin{equation*}
\det(\calT,\calP)\to\det(\calE,\calP).
\end{equation*}
If in addition $\calT$ is small we deduce that $K_i(\calT)\cong
K_i(\calE)$ for $i=0, 1$. The proof that the functor
$\det(\calT,\calP)\to\det(\calE,\calP)$ is an equivalence of
categories contains several intermediate steps of independent
interest. For example we will show that any bounded cohomological
functor on a triangulated category induces a functor between certain
categories of determinant functors.

We should mention two topics which are not addressed in this paper.
Firstly we do not give any applications of determinant functors on
triangulated categories. One application which will be discussed
elsewhere is a new construction of certain Euler characteristics in
relative algebraic $K_0$-groups (which are defined in
\cite{BurnsFlach01}, \cite{Burns04} and \cite{BreuningBurns05})
using universal determinant functors on triangulated categories of
perfect complexes. Secondly many basic questions regarding the
groups $K_0$ and $K_1$ of a triangulated category remain open. In
particular we do not know how the group $K_1$ is related to the
$K_1$-groups of a triangulated category defined by Neeman
\cite{Neeman05}.

Finally we would like to remark that results in the spirit of this
paper have also been shown for Waldhausen categories instead of
triangulated categories. Further details on this can be found in the
paper by Muro and Tonks \cite{MuroTonks06} and in forthcoming work
by Witte.

This paper is organized as follows. In \S
\ref{section_preliminaries} we recall the definitions of Picard
categories, monoidal functors and determinant functors on exact
categories, and we fix our notation for triangulated categories. The
definition and basic properties of determinant functors on
triangulated categories are given in \S \ref{section_definition}.
Next in \S \ref{section_universal} we show the existence of a
universal determinant functor for every small triangulated category
and we define the groups $K_i(\calT)$ for $i=0$ and $1$. Finally in
\S \ref{section_t_structure} we prove the comparison result for
determinant functors on a triangulated category with a
non-degenerate bounded t-structure and determinant functors on its
heart.

I would like to thank David Burns and Matthias Flach for helpful
discussions.

\section{Preliminaries}
\label{section_preliminaries}

In this section we collect some necessary background material. We
first summarize the relevant facts about Picard categories and
monoidal functors between Picard categories. We then recall the
definition of determinant functors on exact categories. Finally we
fix the terminology that we use for triangulated categories.

\subsection{Picard categories}
\label{subsection_Picard}

An AC tensor category is a category $\calP$ together with a
bifunctor $\otimes: \calP\times\calP\to\calP$ and compatible
associativity and commutativity constraints $\varphi$ and $\psi$
respectively, i.e.\ $\varphi_{X,Y,Z}: X\otimes (Y\otimes
Z)\xrightarrow{\cong} (X\otimes Y)\otimes Z$ and $\psi_{X,Y}:
X\otimes Y\xrightarrow{\cong} Y\otimes X$ are natural isomorphisms
satisfying $\psi_{Y,X}\circ\psi_{X,Y}=\id_{X\otimes Y}$, the
pentagonal axiom \cite[Chap.\ I, 1.1.1]{Saavedra72} and the
hexagonal axiom \cite[Chap.\ I, 2.1.1]{Saavedra72}. We will write
$\calP$ for an AC tensor category $(\calP,\otimes,\varphi,\psi)$ if
$\otimes$, $\varphi$ and $\psi$ are clear from the context. The
coherence theorem for AC tensor categories often allows us to omit
explicit reference to $\varphi$ and $\psi$ and to write expressions
involving $\otimes$ without parentheses; see \cite[Chap.\ I,
2.1.2]{Saavedra72} for a precise formulation.

\begin{definition}
An AC tensor category $\calP=(\calP,\otimes,\varphi,\psi)$ is called
a Picard category if it satisfies the following three conditions:
$\calP$ is non-empty, every morphism in $\calP$ is an isomorphism,
and for every object $W$ in $\calP$ the functor $X\mapsto W\otimes
X$ is an autoequivalence of $\calP$.
\end{definition}

We remark that such a category is sometimes called a commutative
Picard category. In this paper we will not consider any
non-commutative Picard categories and we will therefore always omit
the word commutative.

A unit in a Picard category $\calP$ is a triple $(U,\lambda,\rho)$
consisting of an object $U$ and natural isomorphisms $\lambda_X:
X\xrightarrow{\cong} U\otimes X$ and $\rho_X: X\xrightarrow{\cong}
X\otimes U$ satisfying the conditions \cite[Chap.\ I, 1.3.1 and
2.2.1]{Saavedra72} (where the notation $g$ and $d$ is used instead
of $\lambda$ and $\rho$). A morphism of units $(U,\lambda,\rho)\to
(U',\lambda',\rho')$ is a morphism $\alpha: U\to U'$ such that
$(\alpha\otimes\id_X)\circ\lambda_X=\lambda'_X$ and
$(\id_X\otimes\alpha)\circ\rho_X=\rho'_X$ for all objects $X$. A
reduced unit in $\calP$ is a pair $(U,\delta)$ where $U$ is an
object and $\delta: U\xrightarrow{\cong} U\otimes U$ is an
isomorphism. A morphism of reduced units $(U,\delta)\to
(U',\delta')$ is a morphism $\alpha: U\to U'$ such that
$(\alpha\otimes\alpha)\circ\delta = \delta'\circ\alpha$. By
\cite[Chap.\ I, 2.2.5]{Saavedra72} the functor from the category of
units to the category of reduced units which sends
$(U,\lambda,\rho)$ to $(U,\lambda_U)$ (note that $\lambda_U=\rho_U$)
is an isomorphism of categories. We will usually identify a reduced
unit $(U,\delta)$ with the unique unit $(U,\lambda,\rho)$ for which
$\delta=\lambda_U=\rho_U$, and in particular call the pair
$(U,\delta)$ a unit. To simplify the notation we will sometimes
write $\delta$ for the isomorphisms $\lambda_X: X\to U\otimes X$ and
$\rho_X: X\to X\otimes U$.

There exists a natural product of two units. If $(U,\lambda,\rho)$
and $(U',\lambda',\rho')$ are two units then we define their product
$(U,\lambda,\rho)\otimes (U',\lambda',\rho')$ to be the unit
$(U\otimes U',\lambda'',\rho'')$ where $\lambda''_X$ is the
composite $X\xrightarrow{\lambda'_X} U'\otimes
X\xrightarrow{\lambda_{U'}\otimes\id} (U\otimes U')\otimes X$ and
$\rho''_X$ is the composite $X\xrightarrow{\rho_X} X\otimes
U\xrightarrow{\id\otimes\rho'_U} X\otimes (U\otimes U')$. In terms
of reduced units the product has the following description. If
$(U,\delta)$ and $(U',\delta')$ are two (reduced) units in $\calP$
then their product is $(U,\delta)\otimes (U',\delta')=(U\otimes
U',\delta'')$ where $\delta''$ is the composite isomorphism
$U\otimes U'\xrightarrow{\delta\otimes\delta'} (U\otimes U)\otimes
(U'\otimes U')\cong (U\otimes U')\otimes (U\otimes U')$ (the last
isomorphism does not depend on the order of the two copies of $U$
and of $U'$ because $\psi_{U,U}=\id_{U\otimes U}$ and
$\psi_{U',U'}=\id_{U'\otimes U'}$).

Every Picard category has a unit and units are unique up to unique
isomorphism. Therefore computations with units can also be performed
by fixing a unit in $\calP$. If $(1_\calP,\delta)$ with $\delta:
1_\calP\to 1_\calP\otimes 1_\calP$ is a fixed unit in $\calP$ then
the category of units is isomorphic to the category of pairs
$(U,\omega)$ where $U$ is an object in $\calP$ and $\omega: U\to
1_\calP$ is an isomorphism, and where a morphism $(U,\omega)\to
(U',\omega')$ is a morphism $\alpha:U\to U'$ such that
$\omega'\circ\alpha=\omega$. The product of $(U,\omega)$ and
$(U',\omega')$ is given by $(U\otimes U',\omega'')$ where $\omega''$
is the composite $U\otimes U'\xrightarrow{\omega\otimes\omega'}
1_\calP\otimes 1_\calP\xrightarrow{\delta^{-1}} 1_\calP$. Since in
general a Picard category does not have a unique (or a canonical)
unit, we prefer not to work with a fixed unit in this paper.

The homotopy groups $\pi_0(\calP)$ and $\pi_1(\calP)$ of a small
Picard category $\calP$ are defined as follows: $\pi_0(\calP)$ is
the group of isomorphism classes of objects of $\calP$ (with product
induced by $\otimes$) and $\pi_1(\calP)$ is the group of
automorphisms of any object $X$ of $\calP$ (this is well-defined
because for any two objects $X$ and $Y$ there exists a canonical
isomorphism $\Aut_{\calP}(X)\cong\Aut_{\calP}(Y)$). The groups
$\pi_0(\calP)$ and $\pi_1(\calP)$ are abelian. We define a map
$\varepsilon:\pi_0(\calP)\to\pi_1(\calP)$ by sending the isomorphism
class of an object $X$ to the automorphism
$\psi_{X,X}\in\Aut_{\calP}(X\otimes X)=\pi_1(\calP)$. It is not
difficult to show that $\varepsilon$ is a homomorphism.

\subsection{Monoidal functors}
\label{subsection_monoidal}

A monoidal functor from a Picard category
$(\calP,\otimes,\varphi,\psi)$ to a Picard category
$(\calP',\otimes',\varphi',\psi')$ is a pair $(M,c)$ consisting of a
functor $M: \calP\to\calP'$ and a natural isomorphism $c_{X,Y}:
M(X)\otimes' M(Y)\xrightarrow{\cong} M(X\otimes Y)$ which is
compatible with the associativity and commutativity constraints (see
\cite[Chap.\ I, 4.2.1, 4.2.2]{Saavedra72}). We shall write $M$ for
the monoidal functor $(M,c)$ if no confusion can occur. If $(M,c),
(N,d)$ are monoidal functors from $\calP$ to $\calP'$ then a
morphism from $(M,c)$ to $(N,d)$ is a natural transformation
$\kappa: M\to N$ which satisfies $\kappa_{X\otimes Y}\circ
c_{X,Y}=d_{X,Y}\circ(\kappa_X\otimes'\kappa_Y)$. Every morphism of
monoidal functors of Picard categories is automatically an
isomorphism. We write $\Hom^{\otimes}(\calP,\calP')$ for the
category whose objects are monoidal functors from $\calP$ to
$\calP'$ and whose morphisms are morphisms of monoidal functors.
Finally we note that if $(M,c): \calP\to\calP'$ and $(M',c'):
\calP'\to\calP''$ are monoidal functors then there is an obvious way
to define the composition $(M',c')\circ (M,c): \calP\to\calP''$.

A monoidal functor $(M,c): \calP\to\calP'$ of small Picard
categories induces homomorphisms $\pi_i(M):
\pi_i(\calP)\to\pi_i(\calP')$ for $i=0, 1$. It is not difficult to
see that two isomorphic monoidal functors induce the same
homomorphisms on $\pi_i$ for $i=0, 1$.

\begin{lemma}
For a monoidal functor $(M,c): \calP\to\calP'$ of small Picard
categories the following are equivalent.
\begin{itemize}
\item[(i)] $M$ is an equivalence of categories, i.e.\ there exists a
functor $N: \calP'\to\calP$ (not a monoidal functor) such that
$N\circ M$ and $M\circ N$ are isomorphic to the respective identity
functors.
\item[(ii)] $(M,c)$ is an equivalence of Picard categories, i.e.\
there exists a monoidal functor $(N,d): \calP'\to\calP$ such that
$(N,d)\circ (M,c)$ and $(M,c)\circ (N,d)$ are isomorphic to the
respective identity functors (isomorphic as monoidal functors).
\item[(iii)] $(M,c)$ induces isomorphisms on $\pi_0$ and $\pi_1$.
\end{itemize}
\end{lemma}

\begin{proof}
The implication (ii) $\Longrightarrow$ (i) is obvious, and (i)
$\Longrightarrow$ (iii) holds because if $M: \calP\to\calP'$ is an
equivalence of categories then $M$ induces a bijection on the set of
isomorphism classes of $\calP$ and $\calP'$, and an isomorphism
$\Aut_\calP(X)\cong\Aut_{\calP'}(M(X))$ for any object $X$ of
$\calP$.

To show (iii) $\Longrightarrow$ (ii) we construct a monoidal functor
$(N,d): \calP'\to\calP$ as follows. For every object $X'$ of
$\calP'$ we fix an object $N(X')$ in $\calP$ and an isomorphism
$\kappa_{X'}: M(N(X'))\to X'$ in $\calP'$ (such $N(X')$ and
$\kappa_{X'}$ exist because $M$ induces a surjection on isomorphism
classes). Since $M$ induces an injection on isomorphism classes and
a bijection on automorphisms, it follows that for any isomorphism
$\alpha': X'\to Y'$ in $\calP'$ there exists a unique isomorphism
$N(\alpha'): N(X')\to N(Y')$ such that
$\alpha'\circ\kappa_{X'}=\kappa_{Y'}\circ M(N(\alpha'))$.
Furthermore for every pair $X', Y'$ there exists a unique
$d_{X',Y'}: N(X')\otimes N(Y')\to N(X'\otimes' Y')$ such that
$\kappa_{X'}\otimes'\kappa_{Y'}=\kappa_{X'\otimes' Y'}\circ
M(d_{X',Y'})\circ c_{N(X'),N(Y')}$. One easily verifies that $(N,d):
\calP'\to\calP$ is a monoidal functor and that $\kappa: (M,c)\circ
(N,d)\to \id_{\calP'}$ is an isomorphism of monoidal functors.
Moreover it is not difficult to see that there also exists an
isomorphism $(N,d)\circ (M,c)\cong\id_{\calP}$.
\end{proof}

The following lemma describes the automorphisms of a monoidal
functor in terms of homomorphisms between certain homotopy groups.

\begin{lemma}
\label{lemma_monoidal_automorphism} Let $\calP$ and $\calP'$ be
small Picard categories and $(M,c): \calP\to\calP'$ a monoidal
functor. Then there exists a canonical isomorphism
\begin{equation*}
\Aut_{\Hom^\otimes(\calP,\calP')}((M,c))\cong
\Hom(\pi_0(\calP),\pi_1(\calP'))
\end{equation*}
(where the last $\Hom$ is homomorphisms of abelian groups).
\end{lemma}

\begin{proof}
Let $\kappa\in\Aut_{\Hom^\otimes(\calP,\calP')}((M,c))$. Then for
every object $X$ in $\calP$ we have an automorphism
$\kappa_X\in\Aut_{\calP'}(M(X))=\pi_1(\calP')$. We define a
homomorphism $\varphi:\pi_0(\calP)\to\pi_1(\calP')$ by
$\varphi([X])=\kappa_X$, where $[X]$ denotes the class of $X$ in
$\pi_0(\calP)$. This is well-defined as a map of sets because if $X$
and $Y$ are isomorphic in $\calP$ then $\kappa_X=\kappa_Y$ in
$\pi_1(\calP')$ by the naturality of $\kappa$. Furthermore $\varphi$
is a homomorphism because the commutative diagram
\begin{equation*}
\xymatrix@C+0.5cm{ M(X)\otimes' M(Y)
\ar[d]^{\kappa_X\otimes'\kappa_Y}
\ar[r]^-{c_{X,Y}} & M(X\otimes Y) \ar[d]^{\kappa_{X\otimes Y}} \\
M(X)\otimes' M(Y) \ar[r]^-{c_{X,Y}} & M(X\otimes Y) }
\end{equation*}
implies that $\kappa_{X\otimes Y}=\kappa_X\cdot\kappa_Y$ in
$\pi_1(\calP')$.

Conversely, if $\varphi: \pi_0(\calP)\to\pi_1(\calP')$ is a
homomorphism then setting
$\kappa_X=\varphi([X])\in\pi_1(\calP')=\Aut_{\calP'}(M(X))$ for an
object $X$ in $\calP$ defines an automorphism $\kappa$ of $(M,c)$.
Obviously $\kappa\mapsto\varphi$ and $\varphi\mapsto\kappa$ are
mutually inverse isomorphisms.
\end{proof}

\subsection{Determinant functors on exact categories}
\label{subsection_det_exact}

Let $\calE$ be an exact category. We will write the exact sequences
in $\calE$ either as $A\to B\to C$ or as $0\to A\to B\to C\to 0$ and
we will always refer to them as short exact sequences. Let $w$ be a
class of morphisms in $\calE$ which contains all isomorphisms and is
closed under composition. The subcategory of $\calE$ where morphisms
are restricted to $w$ will be denoted by $\calE_w$.

\begin{definition}
\label{definition_det_exact} A determinant functor $f$ on
$(\calE,w)$ with values in a Picard category $\calP$ is a pair
$f=(f_1,f_2)$ consisting of a functor $f_1:\calE_w\to\calP$ and for
every short exact sequence $\Delta: 0\to A\to B\to C\to 0$ in
$\calE$ an isomorphism $f_2(\Delta): f_1(B)\to f_1(A)\otimes
f_1(C)$, such that the following three axioms are satisfied.
\begin{enumerate}
\item[(i)] \textit{Naturality.} For every morphism of short exact sequences
\begin{equation*}
\xymatrix{ \Delta: 0 \ar[r] & A \ar[r] \ar[d]^{a} & B \ar[r]
\ar[d]^{b} & C \ar[r] \ar[d]^{c} & 0 \\
\Delta': 0\ar[r] & A' \ar[r] & B' \ar[r] & C' \ar[r] & 0 }
\end{equation*}
with $a, b, c\in w$, the following diagram is commutative.
\begin{equation*}
\xymatrix@C+1cm{ f_1(B) \ar[r]^-{f_2(\Delta)} \ar[d]^{f_1(b)} &
f_1(A)\otimes f_1(C) \ar[d]^{f_1(a)\otimes f_1(c)} \\
f_1(B') \ar[r]^-{f_2(\Delta')} & f_1(A')\otimes f_1(C') }
\end{equation*}
\item[(ii)] \textit{Associativity.} For every commutative diagram of short exact
sequences
\begin{equation}
\label{diagram_assoc_exact}
\begin{split}
\xymatrix{
A \ar[r] \ar@{=}[d] & B \ar[r] \ar[d] & C' \ar[d] \\
A \ar[r] & C \ar[r] \ar[d] & B' \ar[d] \\
& A' \ar@{=}[r] & A', }
\end{split}
\end{equation}
the following diagram is commutative
\begin{equation*}
\xymatrix@C+1.5cm{ f_1(C) \ar[dd]^{f_2(\Delta_{\mathrm{v1}})}
\ar[r]^-{f_2(\Delta_{\mathrm{h2}})}
& f_1(A)\otimes f_1(B') \ar[d]^{\id\otimes f_2(\Delta_{\mathrm{v2}})} \\
& f_1(A)\otimes (f_1(C')\otimes f_1(A')) \ar[d]^{\varphi_{f_1(A),f_1(C'),f_1(A')}} \\
f_1(B)\otimes f_1(A') \ar[r]^-{f_2(\Delta_{\mathrm{h1}})\otimes\id}
& (f_1(A)\otimes f_1(C'))\otimes f_1(A'), }
\end{equation*}
where $\Delta_{\mathrm{h1}}$ and $\Delta_{\mathrm{h2}}$ (resp.\
$\Delta_{\mathrm{v1}}$ and $\Delta_{\mathrm{v2}}$) are the first and
second horizontal (resp.\ vertical) short exact sequences in diagram
(\ref{diagram_assoc_exact}).
\item[(iii)] \textit{Commutativity.}
For every pair of short exact sequences
\begin{gather*}
\Delta_1: A\to A\oplus B\to B, \\
\Delta_2: B\to A\oplus B\to A
\end{gather*}
where $A\oplus B$ is a direct sum of $A$ and $B$ and the maps are
the canonical inclusions and projections, the following diagram is
commutative.
\begin{equation*}
\xymatrix{
& f_1(A\oplus B) \ar[dl]_{f_2(\Delta_1)} \ar[dr]^{f_2(\Delta_2)} & \\
f_1(A)\otimes f_1(B) \ar[rr]^{\psi_{f_1(A),f_1(B)}} & &
f_1(B)\otimes f_1(A) }
\end{equation*}
\end{enumerate}
\end{definition}

\begin{remark}
Our definition of a determinant functor agrees with the definition
in \cite{Knudsen02} because the compatibility axiom in
\cite[Definition 1.4]{Knudsen02} follows from the associativity
axiom. In \cite{Knudsen02}, $w$ is required to be an SQ-class of
morphisms, however the definition makes sense for any class of
morphisms $w$ as considered above.
\end{remark}

\begin{definition}
Let $f=(f_1,f_2)$ and $g=(g_1,g_2)$ be determinant functors on
$(\calE,w)$ with values in $\calP$. A morphism $f\to g$ is a natural
isomorphism $\lambda: f_1\to g_1$ such that
$(\lambda_A\otimes\lambda_C)\circ f_2(\Delta) =
g_2(\Delta)\circ\lambda_B: f_1(B)\to g_1(A)\otimes g_1(C)$ for every
short exact sequence $\Delta: 0\to A\to B\to C\to 0$ in $\calE$.
\end{definition}

We will often use the notation $f: (\calE,w)\to\calP$ to denote a
determinant functor $f$ on $(\calE,w)$ with values in $\calP$. We
write $\det((\calE,w),\calP)$ for the category whose objects are the
determinant functors $(\calE,w)\to\calP$ and whose morphisms are the
morphisms of determinant functors. If $F:(\calE',w')\to(\calE,w)$ is
an exact functor (i.e.\ $F:\calE'\to\calE$ is an exact functor of
exact categories which satisfies $F(w')\subseteq w$) and
$M:\calP\to\calP'$ is a monoidal functor then we obtain an induced
functor $\det((\calE,w),\calP)\to\det((\calE',w'),\calP')$, which
sends a determinant functor $f: (\calE,w)\to\calP$ to the composite
determinant functor $M\circ f\circ F: (\calE',w')\to\calP'$ that is
defined in the obvious way.

If $w=\iso$ is the class of all isomorphisms in $\calE$ then we will
in general omit $w$ from the notation. So we write
$\det(\calE,\calP)$ instead of $\det((\calE,\iso),\calP)$ and talk
about determinant functors on $\calE$ with values in $\calP$ instead
of determinant functors on $(\calE,\iso)$ with values in $\calP$.

\subsection{Preliminaries about triangulated categories}

Let $\calT$ be a triangulated category, i.e.\ $\calT$ is an additive
category with an additive automorphism $T$ (the translation functor)
and a class of distinguished triangles $A\to B\to C\to TA$ such that
the axioms (TR1)--(TR4) in \cite{Verdier} hold.

Let $A$ and $B$ be objects in $\calT$. It follows from the axioms of
a triangulated category that
\begin{equation}
\label{eqn_direct_sum} A\to A\oplus B\to B\xrightarrow{0} TA
\end{equation}
is a distinguished triangle in $\calT$, where $A\oplus B$ is a
direct sum of $A$ and $B$, and the maps $A\to A\oplus B$ and
$A\oplus B\to B$ are the canonical inclusion and projection
respectively.

By an octahedral diagram in $\calT$ we mean a diagram of the form
\begin{equation}
\label{eqn_octahedral_diagram}
\begin{split}
\xymatrix{ A \ar[r]^a \ar@{=}[d] & B \ar[r]^{b} \ar[d] & C' \ar[r]
\ar[d] & TA \ar@{=}[d] \\
A \ar[r] & C \ar[r] \ar[d] & B' \ar[r]^c \ar[d]^d & TA \\
& A' \ar@{=}[r] \ar[d]^e & A' \ar[d] & \\
& TB \ar[r]^{Tb} & TC' & }
\end{split}
\end{equation}
such that
\begin{enumerate}
\item[(i)] the diagram is commutative,
\item[(ii)] the two horizontal triangles and the two vertical triangles
are distinguished,
\item[(iii)] the two composite morphisms $B'\xrightarrow{c} TA\xrightarrow{Ta} TB$ and
$B'\xrightarrow{d} A'\xrightarrow{e} TB$ are equal.
\end{enumerate}
We remark that octahedral diagrams are the diagrams which appear in
the octahedral axiom (TR4) in \cite{Verdier}.

An exact functor $(F,\tau): \calS\to\calT$ of triangulated
categories is a pair consisting of an additive functor
$F:\calS\to\calT$ and a natural isomorphism $\tau: F\circ T_\calS\to
T_\calT\circ F$ such that if $A\xrightarrow{a} B\xrightarrow{b}
C\xrightarrow{c} T_{\calS}A$ is a distinguished triangle in $\calS$
then $FA\xrightarrow{Fa} FB\xrightarrow{Fb}
FC\xrightarrow{\tau_A\circ Fc} T_{\calT}FA$ is a distinguished
triangle in $\calT$. An exact functor sends a distinguished triangle
of the form (\ref{eqn_direct_sum}) in $\calS$ to a distinguished
triangle of the same form in $\calT$. Furthermore it sends an
octahedral diagram in $\calS$ to an octahedral diagram in $\calT$.
For the definition of a morphism of exact functors we refer to
\cite[Chap.\ 2, \S 2, no.\ 1]{Verdier}.

\section{Definition and basic properties}
\label{section_definition}

In this section we define determinant functors on triangulated
categories and prove various basic properties for them.

\subsection{Definition of determinant functors}

Let $\calT$ be a triangulated category. We write $\iso$ for the
class of isomorphisms in $\calT$ and $\calT_{\iso}$ for the
subcategory of $\calT$ where morphisms are restricted to
isomorphisms.

\begin{definition}
A determinant functor $f$ on $\calT$ with values in a Picard
category $\calP$ is a pair $f=(f_1,f_2)$ consisting of a functor
$f_1: \calT_{\iso}\to\calP$ and for every distinguished triangle
$\Delta: A\to B\to C\to TA$ in $\calT$ an isomorphism $f_2(\Delta):
f_1(B)\to f_1(A)\otimes f_1(C)$, such that the following three
axioms are satisfied.
\begin{enumerate}
\item[(i)] \textit{Naturality.} For every isomorphism of distinguished
triangles
\begin{equation}
\label{diagram_iso_dist} \begin{split} \xymatrix{ \Delta{\,}:{\,}
A{\,} \ar[r] \ar@<2.3ex>[d]^{a} & B \ar[r] \ar[d]^{b} & C
\ar[r] \ar[d]^{c} & TA \ar[d]^{Ta} \\
\Delta':{\,} A' \ar[r] & B' \ar[r] & C' \ar[r] & TA', }
\end{split}
\end{equation}
the following diagram is commutative.
\begin{equation*}
\xymatrix@C+1cm{ f_1(B) \ar[r]^-{f_2(\Delta)} \ar[d]^{f_1(b)} &
f_1(A)\otimes f_1(C) \ar[d]^{f_1(a)\otimes f_1(c)} \\
f_1(B') \ar[r]^-{f_2(\Delta')} & f_1(A')\otimes f_1(C') }
\end{equation*}
\item[(ii)] \textit{Associativity.} For every octahedral diagram as in
(\ref{eqn_octahedral_diagram}), the following diagram is commutative
\begin{equation*}
\xymatrix@C+1.5cm{ f_1(C) \ar[dd]^{f_2(\Delta_{\mathrm{v1}})}
\ar[r]^-{f_2(\Delta_{\mathrm{h2}})}
& f_1(A)\otimes f_1(B') \ar[d]^{\id\otimes f_2(\Delta_{\mathrm{v2}})} \\
& f_1(A)\otimes (f_1(C')\otimes f_1(A')) \ar[d]^{\varphi_{f_1(A),f_1(C'),f_1(A')}} \\
f_1(B)\otimes f_1(A') \ar[r]^-{f_2(\Delta_{\mathrm{h1}})\otimes\id}
& (f_1(A)\otimes f_1(C'))\otimes f_1(A'), }
\end{equation*}
where $\Delta_{\mathrm{h1}}$ and $\Delta_{\mathrm{h2}}$ (resp.\
$\Delta_{\mathrm{v1}}$ and $\Delta_{\mathrm{v2}}$) are the first and
second horizontal (resp.\ vertical) distinguished triangles in
diagram (\ref{eqn_octahedral_diagram}).
\item[(iii)] \textit{Commutativity.}
For every pair of distinguished triangles
\begin{equation}
\label{eqn_split_dist}
\begin{split}
\Delta_1: & A\to A\oplus B\to B\xrightarrow{0} TA, \\
\Delta_2: & B\to A\oplus B\to A\xrightarrow{0} TB
\end{split}
\end{equation}
where $A\oplus B$ is a direct sum of $A$ and $B$ and the maps are
the canonical inclusions and projections, the following diagram is
commutative.
\begin{equation*}
\xymatrix{
& f_1(A\oplus B) \ar[dl]_{f_2(\Delta_1)} \ar[dr]^{f_2(\Delta_2)} & \\
f_1(A)\otimes f_1(B) \ar[rr]^{\psi_{f_1(A),f_1(B)}} & &
f_1(B)\otimes f_1(A) }
\end{equation*}
\end{enumerate}
\end{definition}

\begin{definition}
Let $f=(f_1,f_2)$ and $g=(g_1,g_2)$ be determinant functors on
$\calT$ with values in $\calP$. A morphism $f\to g$ is a natural
isomorphism $\lambda: f_1\to g_1$ such that
$(\lambda_A\otimes\lambda_C)\circ f_2(\Delta) =
g_2(\Delta)\circ\lambda_B: f_1(B)\to g_1(A)\otimes g_1(C)$ for every
distinguished triangle $\Delta: A\to B\to C\to TA$ in $\calT$.
\end{definition}

We will often use the notation $f: \calT\to\calP$ to denote a
determinant functor $f$ on $\calT$ with values in $\calP$. We write
$\det(\calT,\calP)$ for the category whose objects are the
determinant functors $\calT\to\calP$ and whose morphisms are the
morphisms of determinant functors. If $F:\calT'\to\calT$ is an exact
functor of triangulated categories and $M:\calP\to\calP'$ is a
monoidal functor then we obtain an induced functor
$\det(\calT,\calP)\to\det(\calT',\calP')$, which sends a determinant
functor $f: \calT\to\calP$ to the composite determinant functor
$M\circ f\circ F: \calT'\to\calP'$ that is defined in the obvious
way.

To simplify the notation we will often write $[\cdot]$ for a
determinant functor $f=(f_1,f_2)$, i.e.\ if $A$ (resp.\ $a$) is an
object (resp.\ isomorphism) in $\calT$ then $[A]=f_1(A)$ (resp.\
$[a]=f_1(a)$), and if $\Delta$ is a distinguished triangle in
$\calT$ then $[\Delta]=f_2(\Delta)$. This will not cause any
confusion.

\subsection{Basic properties of determinant functors}
\label{subsection_basic_properties}

Let $[\cdot]$ be a determinant functor on a triangulated category
$\calT$ with values in a Picard category $\calP$.

\begin{lemma}
\label{lemma_zero} If $0$ is a zero object in $\calT$ then $[0]$ has
canonically the structure of a unit in $\calP$, i.e.\ there exists a
canonical isomorphism $\delta_0: [0]\to [0]\otimes [0]$. If $a: 0\to
0'$ is the unique isomorphism of zero objects $0$ and $0'$ then
$[a]: [0]\to [0']$ is an isomorphism of units $([0],\delta_0)\to
([0'],\delta_{0'})$.
\end{lemma}

\begin{proof}
Applying $[\cdot]$ to the distinguished triangle $0\to 0\to 0\to T0$
gives a canonical isomorphism $\delta_0: [0]\to [0]\otimes [0]$ and
therefore the structure of a unit on $[0]$. The isomorphism $a: 0\to
0'$ induces an isomorphism from the distinguished triangle $0\to
0\to 0\to T0$ to $0'\to 0'\to 0'\to T0'$, and therefore $[a]$
commutes with the unit structures on $[0]$ and $[0']$.
\end{proof}

We remark that the unit structure $\delta_0$ on $[0]$ depends on the
zero object $0$ in $\calT$ and not only on the object $[0]$ in
$\calP$. Nevertheless to simplify the notation we will talk about
the unit $[0]$ when we mean the unit $([0],\delta_0)$.

\begin{lemma}
\label{lemma_compatibility} Let $a: A\to B$ be an isomorphism and
$0$ a zero object in $\calT$. Let $\lambda_X: X\to [0]\otimes X$ and
$\rho_{X}: X\to X\otimes [0]$ be the natural isomorphisms associated
to the unit $[0]$.
\begin{itemize}
\item[(i)] Let $\Delta$ be the distinguished triangle $A\xrightarrow{a}
B\to 0 \to TA$. Then $[a]=[\Delta]^{-1}\circ \rho_{[A]}: [A]\to
[B]$.
\item[(ii)] Let $\Delta$ be the distinguished triangle $0\to
A\xrightarrow{a} B \to T0$. Then $[a]=\lambda^{-1}_{[B]}\circ
[\Delta]: [A]\to [B]$.
\end{itemize}
\end{lemma}

\begin{proof}
We prove part (ii), the proof of part (i) is similar and left to the
reader. Recall that by definition of the unit structure on $[0]$ the
isomorphism $\delta_0=\lambda_{[0]}=\rho_{[0]}: [0]\to [0]\otimes
[0]$ is induced by the distinguished triangle $0\to 0\to 0\to T0$.
Now consider the distinguished triangle $\Delta': 0\to
B\xrightarrow{\id} B\to T0$. From the octahedral diagram
\begin{equation*}
\xymatrix{
0 \ar[r] \ar@{=}[d] & 0 \ar[r] \ar[d] & 0 \ar[r] \ar[d] & T0 \ar@{=}[d] \\
0 \ar[r] & B \ar[r]^{\id} \ar[d]^{\id} & B \ar[r] \ar[d]^{\id} & T0 \\
& B \ar@{=}[r] \ar[d] & B \ar[d] & \\
& T0 \ar[r] & T0 & }
\end{equation*}
we obtain $(\delta_0\otimes\id_{[B]})\circ [\Delta']=
\varphi_{[0],[0],[B]}\circ (\id_{[0]}\otimes [\Delta'])\circ
[\Delta']$. On the other hand by the axioms for a unit
$\delta_0\otimes\id_{[B]}=\varphi_{[0],[0],[B]}\circ
(\id_{[0]}\otimes \lambda_{[B]})$. The last two equations imply
$[\Delta']=\lambda_{[B]}$. The isomorphism of distinguished
triangles
\begin{equation*}
\xymatrix{ 0 \ar[r] \ar[d] & A \ar[r]^{a} \ar[d]^{a} & B \ar[r]
\ar[d]^{\id}
& T0 \ar[d] \\
0 \ar[r] & B \ar[r]^{\id} & B \ar[r] & T0 }
\end{equation*}
shows that $[\Delta']\circ[a]=[\Delta]$. Therefore
$[a]=\lambda^{-1}_{[B]}\circ[\Delta]$ as required.
\end{proof}

\begin{lemma}
\label{lemma_minus} Let $A$ be an object of $\calT$ and let $-\id_A:
A\to A$ be the negative of the identity on $A$. Then
$[-\id_A]=\varepsilon([A])$ in $\pi_1(\calP)$.
\end{lemma}

\begin{proof}
The proof is essentially the same as in the case of exact categories
(cf.\ \cite[\S 4.9]{Deligne87} or \cite[Prop.\ 1.5(c)]{Knudsen02}).

Recall that by definition $\varepsilon([A])=\psi_{[A],[A]}$ in
$\pi_1(\calP)$. If $d: A\oplus A\to A\oplus A$ is the isomorphism
$d(a,a')=(a',a)$ then it follows from the naturality and
commutativity axioms that $\psi_{[A],[A]}=[d]$ in $\pi_1(\calP)$,
and from the isomorphism of distinguished triangles
\begin{equation*}
\xymatrix@C+1cm{ A \ar[r]^-{(\id,\id)} \ar[d]^{\id} & A\oplus A
\ar[r]^-{\langle\id,-\id\rangle}
\ar[d]^{d} & A \ar[r]^{0} \ar[d]^{-\id_A} & TA \ar[d]^{\id} \\
A \ar[r]^-{(\id,\id)} & A\oplus A \ar[r]^-{\langle\id,-\id\rangle} &
A \ar[r]^{0} & TA }
\end{equation*}
we deduce that $[d]=[-\id_A]$ in $\pi_1(\calP)$.
\end{proof}

\begin{lemma}
\label{lemma_mu} If $A$ is an object in $\calT$ then $[A]\otimes
[TA]$ has canonically the structure of a unit in $\calP$, i.e.\
there exists a canonical isomorphism $\mu_A: [A]\otimes [TA]\to
([A]\otimes [TA])\otimes ([A]\otimes [TA])$. This unit structure has
the following properties.
\begin{enumerate}
\item[(i)] If $a: A\to B$ is an isomorphism then $[a]\otimes [Ta]$ is an
isomorphism of units $([A]\otimes [TA],\mu_A)\to ([B]\otimes
[TB],\mu_B)$.
\item[(ii)] If $\Delta: A\xrightarrow{a} B\xrightarrow{b} C\xrightarrow{c}
TA$ is a distinguished triangle and $\Delta'$ denotes the
distinguished triangle $TA \xrightarrow{Ta} TB \xrightarrow{Tb} TC
\xrightarrow{-Tc} T^2 A$, then the map
\begin{equation*}
\begin{split}
[B]\otimes [TB] & \xrightarrow{[\Delta]\otimes [\Delta']}
([A]\otimes [C])\otimes ([TA]\otimes [TC]) \\
& \xrightarrow{\;\cong\;} ([A]\otimes [TA])\otimes ([C]\otimes [TC])
\end{split}
\end{equation*}
is an isomorphism of units $([B]\otimes [TB],\mu_B)\to ([A]\otimes
[TA],\mu_A)\otimes ([C]\otimes [TC],\mu_C)$.
\end{enumerate}
\end{lemma}

To simplify the notation we will always write $\mu_A$ for the
natural isomorphisms $X\to ([A]\otimes [TA])\otimes X$ and $X\to
X\otimes ([A]\otimes [TA])$ associated to the unit $([A]\otimes
[TA],\mu_A)$. The fact that $[A]\otimes [TA]$ has canonically the
structure of a unit can also be expressed by saying that $[TA]$ is
canonically a right inverse of $[A]$.

\begin{proof}[Proof of Lemma \ref{lemma_mu}]
From the distinguished triangle $\Delta_A: A\to 0\to
TA\xrightarrow{\id} TA$ we obtain an isomorphism $[\Delta_A]:
[0]\cong [A]\otimes [TA]$, and by Lemma \ref{lemma_zero} the object
$[0]$ has a unit structure $\delta_0$. Hence $[A]\otimes [TA]$ has
an induced unit structure $\mu_A=([\Delta_A]\otimes
[\Delta_A])\circ\delta_0\circ [\Delta_A]^{-1}$. Clearly $\mu_A$ is
independent of the choice of zero object $0$.

Property (i) holds because an isomorphism $a: A\to B$ induces an
isomorphism from the distinguished triangle $\Delta_A$ to
$\Delta_B$.

For property (ii) it suffices to show the commutativity of the
following diagram.
\begin{equation*}
\xymatrix@C+1.5cm{ [B]\otimes [TB] \ar[r]^-{[\Delta]\otimes
[\Delta']} \ar[d]^{[\Delta_B]^{-1}} & ([A]\otimes [C])\otimes
([TA]\otimes [TC]) \ar[dd]^{\cong} \\
[0] \ar[d]^{\delta_0} & \\
[0]\otimes [0] \ar[r]^-{[\Delta_A]\otimes [\Delta_C]} & ([A]\otimes
[TA])\otimes ([C]\otimes [TC]) }
\end{equation*}
This diagram is the outer square in diagram (\ref{diagram_mu})
below. In (\ref{diagram_mu}), $\tilde{\Delta}$ denotes the
distinguished triangle
$C\xrightarrow{c}TA\xrightarrow{Ta}TB\xrightarrow{Tb}TC$ and
$\rho_X: X\to X\otimes [0]$ is the natural isomorphism associated to
the unit $[0]$.
\begin{equation}
\label{diagram_mu}
\begin{split}
\xymatrix@C+1.5cm{ {\scriptstyle [B]\otimes[TB]}
\ar[r]^-{[\Delta]\otimes\id} \ar[d]^{[\Delta_B]^{-1}} &
{\scriptstyle [A]\otimes[C]\otimes[TB]}
\ar[r]^-{\id\otimes\id\otimes[\Delta']} \ar[d]^{\id\otimes
[\tilde{\Delta}]^{-1}} & {\scriptstyle
([A]\otimes[C])\otimes ([TA]\otimes[TC])} \ar[d]^{\cong} \\
{\scriptstyle [0]} \ar[r]^-{[\Delta_A]} \ar[d]^{\rho_{[0]}=\delta_0}
& {\scriptstyle [A]\otimes[TA]} \ar[r]^-{\mu_C}
\ar[d]^{\rho_{[A]\otimes [TA]}} &
{\scriptstyle ([A]\otimes [TA])\otimes ([C]\otimes [TC])} \ar@{=}[d] \\
{\scriptstyle [0]\otimes [0]} \ar[r]^-{[\Delta_A]\otimes\id} &
{\scriptstyle ([A]\otimes [TA])\otimes [0]}
\ar[r]^-{\id\otimes\id\otimes [\Delta_C]} & {\scriptstyle
([A]\otimes [TA])\otimes ([C]\otimes [TC])} }
\end{split}
\end{equation}
The commutativity of the top left hand square in (\ref{diagram_mu})
follows from the octahedral diagram
\begin{equation*}
\xymatrix{ A \ar[r]^{a} \ar@{=}[d] & B \ar[r]^{b} \ar[d] &
C \ar[r]^{c} \ar[d]^{c} & TA \ar@{=}[d] \\
A \ar[r] & 0 \ar[r] \ar[d] & TA \ar[r]^{\id} \ar[d]^{Ta} & TA \\
& TB \ar@{=}[r] \ar[d]^{\id} & TB \ar[d]^{Tb} & \\
& TB \ar[r]^{Tb} & TC. & }
\end{equation*}
The commutativity of the top right hand square in (\ref{diagram_mu})
follows from Lemma \ref{lemma_rotation} below applied to
$\tilde{\Delta}$ and its rotation $\Delta'$ (note that the proof of
Lemma \ref{lemma_rotation} does not use part (ii) of Lemma
\ref{lemma_mu}). The bottom left hand square in (\ref{diagram_mu})
is commutative by the naturality of $\rho_X$, and the bottom right
hand square by definition of the unit structure on $[C]\otimes
[TC]$.
\end{proof}

\begin{remark}
Let $0$ be a zero object in $\calT$. Lemma \ref{lemma_zero} gives a
unit structure on $[0]$ and on $[T0]$ and we can therefore consider
the product unit structure on $[0]\otimes [T0]$. On the other hand
Lemma \ref{lemma_mu} also gives a unit structure on $[0]\otimes
[T0]$. It is not difficult to see that these two unit structures
coincide.
\end{remark}

\begin{lemma}
\label{lemma_rotation} Let $\Delta: A\xrightarrow{a}
B\xrightarrow{b} C\xrightarrow{c} TA$ be a distinguished triangle
and $\Delta': B\xrightarrow{b} C\xrightarrow{c} TA\xrightarrow{-Ta}
TB$ the rotated triangle. Then the composite isomorphism
\begin{equation*}
\begin{split}
[B] & \xrightarrow{\,\mu_A\,} ([A]\otimes [TA])\otimes [B] \\
& \xrightarrow{\;\cong\;} [A]\otimes ([B]\otimes [TA]) \\
& \xrightarrow{\id\otimes [\Delta']^{-1}} [A]\otimes [C]
\end{split}
\end{equation*}
is equal to $[\Delta]: [B]\to [A]\otimes [C]$.
\end{lemma}

\begin{proof}
Let $\alpha: [B]\to [A]\otimes [C]$ be the displayed isomorphism in
the statement of the lemma. We must show that
$[\Delta]\circ\alpha^{-1}=\id_{[A]\otimes [C]}$. Consider the
automorphism $\beta$ of $[C]$ which is defined as the composite
\begin{equation*}
[C] \xrightarrow{[\Delta']} [B]\otimes [TA]
\xrightarrow{[\Delta]\otimes\id} ([A]\otimes [C])\otimes [TA] \cong
([A]\otimes [TA])\otimes [C] \xrightarrow{\mu_A^{-1}} [C].
\end{equation*}
It is not hard to see that $\varepsilon([A])\circ
(\id_{[A]}\otimes\beta)=[\Delta]\circ\alpha^{-1}$, therefore it
suffices to show that $\beta=\varepsilon([A])$ in $\pi_1(\calP)$.

Consider the distinguished triangles
\begin{gather*}
\Delta_1: A \xrightarrow{0} C \xrightarrow{(\id,0)} C\oplus TA
\xrightarrow{\langle 0,\id\rangle} TA, \\
\Delta_2: C \xrightarrow{(\id,c)} C\oplus TA \xrightarrow{\langle
c,-\id\rangle} TA \xrightarrow{0} TC, \\
\Delta_3: C \xrightarrow{(\id,0)} C\oplus TA \xrightarrow{\langle
0,\id\rangle} TA \xrightarrow{0} TC, \\
\Delta_4: TA \xrightarrow{(0,\id)} C\oplus TA
\xrightarrow{\langle\id, 0\rangle} C \xrightarrow{0} T^2A.
\end{gather*}
The triangles $\Delta_2$ and $\Delta_3$ are isomorphic by the
following diagram
\begin{equation}
\label{diagram_rotation_2}
\begin{split}
\xymatrix@C+0.5cm{ \Delta_3:{\,} C \ar[r]^-{(\id,0)}
\ar@<2.6ex>[d]^{\id} & C\oplus TA \ar[r]^-{\langle 0,\id\rangle}
\ar[d]^{d} & TA \ar[r]^{0} \ar[d]^{\id} & TC \ar[d]^{\id} \\
\Delta_2:{\,} C \ar[r]^-{(\id,c)} & C\oplus TA \ar[r]^-{\langle
c,-\id\rangle} & TA \ar[r]^{0} & TC, }
\end{split}
\end{equation}
where $d: C\oplus TA\to C\oplus TA$ is the map $d(u,v)=(u,c(u)-v)$.
We have $[d]=\varepsilon([A])$ in $\pi_1(\calP)$ because of the
commutative diagram of distinguished triangles
\begin{equation*}
\xymatrix@C+0.5cm{ TA \ar[r]^-{(0,\id)} \ar[d]^{-\id} & C\oplus TA
\ar[r]^-{\langle\id,0\rangle} \ar[d]^{d} &
C \ar[r]^{0} \ar[d]^{\id} & T^2 A \ar[d]^{-\id} \\
TA \ar[r]^-{(0,\id)} & C\oplus TA \ar[r]^-{\langle\id,0\rangle} & C
\ar[r]^{0} & T^2 A, }
\end{equation*}
Lemma \ref{lemma_minus}, and the fact that
$\varepsilon([A])=\varepsilon([TA])$.

From the octahedral diagram
\begin{equation*}
\begin{split}
\xymatrix@C+0.5cm{ A \ar[r]^{a} \ar@{=}[d] & B \ar[r]^{b} \ar[d]^{b}
& C \ar[r]^{c} \ar[d]^{(\id,c)} & TA
\ar@{=}[d] \\
A \ar[r]^{0} & C \ar[r]^-{(\id,0)} \ar[d]^{c} & C\oplus TA
\ar[r]^-{\langle 0,\id\rangle} \ar[d]^{\langle c,-\id\rangle} & TA \\
& TA \ar@{=}[r] \ar[d]^{-Ta} & TA \ar[d]^{0} &
\\
& TB \ar[r]^{Tb} & TC & }
\end{split}
\end{equation*}
we obtain
\begin{equation*}
([\Delta]\otimes\id_{[TA]})\circ [\Delta'] = \varphi\circ
(\id_{[A]}\otimes [\Delta_2])\circ [\Delta_1].
\end{equation*}
Furthermore by (\ref{diagram_rotation_2}) we have
$[\Delta_2]=[\Delta_3]\circ \varepsilon([A])$, and from the
commutativity axiom it follows that $[\Delta_4]=\psi\circ
[\Delta_3]$. To show that $\beta=\varepsilon([A])$ it therefore
suffices to show that
\begin{equation*}
[C] \xrightarrow{[\Delta_1]} [A]\otimes [C\oplus TA]
\xrightarrow{\id\otimes [\Delta_4]} [A]\otimes ([TA]\otimes
[C])\cong ([A]\otimes [TA])\otimes [C] \xrightarrow{\mu_A^{-1}} [C]
\end{equation*}
is the identity map. This follows easily from the octahedral diagram
\begin{equation*}
\xymatrix@C+0.5cm{ A\ar[r] \ar@{=}[d] & 0 \ar[r] \ar[d] & TA
\ar[r]^{\id} \ar[d]^{(0,\id)} & TA \ar@{=}[d] \\
A\ar[r]^{0} & C \ar[r]^-{(\id,0)} \ar[d]^{\id} & C\oplus TA
\ar[r]^-{\langle 0,\id\rangle} \ar[d]^{\langle\id, 0\rangle} & TA \\
& C \ar@{=}[r] \ar[d] & C \ar[d]^{0} & \\
& T0 \ar[r] & T^2 A & }
\end{equation*}
and Lemma \ref{lemma_compatibility}(ii) applied to the distinguished
triangle $0\to C\xrightarrow{\id} C\to T0$.
\end{proof}

\begin{lemma}
Every commutative square
\begin{equation*}
\xymatrix{
A' \ar[r] \ar[d] & B' \ar[d] \\
A \ar[r] & B }
\end{equation*}
in $\calT$ can be completed to a $9$-term diagram
\begin{equation}
\label{diagram_9_term}
\begin{split}
\xymatrix{ A' \ar[r] \ar[d] & B' \ar[r] \ar[d] & C' \ar[r] \ar[d] &
TA' \ar@{-->}[d]
\\
A \ar[r] \ar[d] & B \ar[r] \ar[d] & C \ar[r] \ar[d] & TA \ar@{-->}[d] \\
A'' \ar[r] \ar[d] & B'' \ar[r] \ar[d] & C'' \ar[r] \ar[d]
\ar@{}[dr]|{-} & TA''
\ar@{-->}[d] \\
TA' \ar@{-->}[r] & TB' \ar@{-->}[r] & TC' \ar@{-->}[r] & T^2A' }
\end{split}
\end{equation}
(i.e.\ the dashed arrows are the translates of the continuous
arrows, all squares are commutative except for the bottom right hand
square which is anti-commutative, the first three rows and the first
three columns are distinguished triangles) such that for every
determinant functor $[\cdot]$ on $\calT$ with values in a Picard
category $\calP$ the induced diagram in $\calP$
\begin{equation}
\label{diagram_det_9_term}
\begin{split}
\xymatrix{
[B] \ar[r] \ar[dd] & [A]\otimes [C] \ar[d] \\
& ([A']\otimes [A''])\otimes ([C']\otimes [C'']) \ar[d]^{\cong} \\
[B']\otimes [B''] \ar[r] & ([A']\otimes [C'])\otimes ([A'']\otimes
[C'']) }
\end{split}
\end{equation}
is commutative.
\end{lemma}

\begin{proof}
In the proof of \cite[Prop.\ 1.1.11]{BBD} a $9$-term diagram of the
form (\ref{diagram_9_term}) is constructed by combining three
octahedral diagrams. The commutativity of diagram
(\ref{diagram_det_9_term}) then follows by applying the
associativity axiom to these octahedral diagrams (using Lemma
\ref{lemma_rotation} to deal with the rotated triangles in the
octahedral diagrams). We leave the details to the reader.
\end{proof}

\section{Universal determinant functors and $K$-groups}
\label{section_universal}

In this section we show that every small triangulated category
$\calT$ has a universal determinant functor. The groups $K_i(\calT)$
for $i=0,1$ are then defined as invariants of this universal
determinant functor. The results in this section are motivated by
Deligne's theory of universal determinant functors and virtual
objects for exact categories which we recall first.

\subsection{Universal determinant functors on exact categories}
\label{subsection_universal_exact}

Let $\calE$ be an exact category and $w$ a class of morphisms in
$\calE$ as in \S \ref{subsection_det_exact}.

\begin{definition}
A universal determinant functor on $(\calE,w)$ is a Picard category
$\calV$ together with a determinant functor $f: (\calE,w)\to\calV$,
such that for every Picard category $\calP$ the functor
\begin{equation*}
\Hom^\otimes(\calV,\calP)\to\det\big((\calE,w),\calP\big)
\end{equation*}
which sends $M$ to $M\circ f$ is an equivalence of categories.
\end{definition}

If they exist, universal determinant functors are unique up to
non-unique isomorphism, more precisely if $f: (\calE,w)\to\calV$ and
$f': (\calE,w)\to\calV'$ are universal determinant functors on
$(\calE,w)$ then there exists a monoidal functor $M: \calV\to\calV'$
(which is determined up to a non-unique isomorphism) such that the
determinant functors $M\circ f$ and $f'$ are (non-uniquely)
isomorphic. Furthermore this functor $M: \calV\to\calV'$ is an
equivalence of Picard categories.

Deligne showed in \cite[\S 4]{Deligne87} that there exists a
universal determinant functor on every small exact category $\calE$
(i.e.\ on $(\calE,\iso)$). In fact he claims that the existence
follows from ``un argument standard'', however he does not indicate
what this standard argument is. On the other hand Deligne also gives
a topological construction of a universal determinant functor: using
Quillen's $Q$-construction he constructs a Picard category
$V(\calE)$ (``cat\'egorie des objets virtuels'') and a determinant
functor $[\cdot]: \calE\to V(\calE)$, and he then sketches a proof
that it is universal. This construction shows that there exist
canonical isomorphisms $\pi_i(V(\calE))\cong K_i(\calE)$ for $i=0,
1$, where $K_i(\calE)$ are the Quillen $K$-groups of $\calE$.

In the case where $\calE$ is the exact category of finitely
generated projective modules over a ring, a different construction
of a universal determinant functor on $\calE$ was sketched by Fukaya
and Kato in \cite[\S 1.2]{FukayaKato05}. However some details in
their construction are omitted, and it is not obvious that one
indeed obtains a determinant functor in the sense of Definition
\ref{definition_det_exact}.

The proofs for the existence of a universal determinant functor by
Deligne and by Fukaya and Kato are not very detailed. In \S
\ref{subsection_universal_triang} we will prove in detail the
existence of a universal determinant functor on a triangulated
category; it seems likely that a variant of that argument would also
work for exact categories.

\subsection{Universal determinant functors on triangulated categories}
\label{subsection_universal_triang}

Let $\calT$ be a triangulated category.

\begin{definition}
A universal determinant functor on $\calT$ is a Picard category
$\calV$ together with a determinant functor $f: \calT\to\calV$, such
that for every Picard category $\calP$ the functor
\begin{equation*}
\Hom^\otimes(\calV,\calP)\to\det(\calT,\calP)
\end{equation*}
which sends $M$ to $M\circ f$ is an equivalence of categories.
\end{definition}

If they exist, universal determinant functors on triangulated
categories are unique up to non-unique isomorphism, i.e.\ they
satisfy the same uniqueness property as the universal determinant
functors on exact categories in \S \ref{subsection_universal_exact}.

\begin{theorem}
\label{theorem_existence} Let $\calT$ be a small triangulated
category. Then there exists a universal determinant functor $f:
\calT\to\calV$.
\end{theorem}

\begin{proof}
We proof the theorem by constructing $\calV$ with generators and
relations. Since $\calV$ is a Picard category and therefore carries
a lot of structure, the construction and the verification of the
universal property are quite long and will occupy the remainder of
\S \ref{subsection_universal_triang}.

\textbf{Outline of the construction.} Objects and morphisms of
$\calV$ are defined using words over certain alphabets. An alphabet
$\AA$ is a set and a word over $\AA$ is a finite sequence of
elements of $\AA$. We denote the set of words over $\AA$ by $\AA^*$.
If $A\in\AA$ then we also write $A$ for the word ``$A$'' consisting
of the one element $A$. If $X,Y\in\AA^*$ then $XY$ denotes the
concatenation of $X$ and $Y$. For example in the construction of the
objects of $\calV$ below, $\langle X\boxtimes Y\rangle$ denotes the
word which is the concatenation of the words ``$\langle$'', $X$,
``$\boxtimes$'', $Y$ and ``$\rangle$''.

The construction of the Picard category $\calV$ proceeds in several
steps. We first construct a set $\objects(\calV)$ which will be the
set of objects of $\calV$. For every pair $X,Y\in\objects(\calV)$ we
then construct a set $\Hom_\calV(X,Y)$ which will be the set of
morphisms from $X$ to $Y$ in $\calV$. We define the composition of
morphisms $\circ: \Hom_\calV(Y,Z)\times\Hom_\calV(X,Y) \to
\Hom_\calV(X,Z)$ so that we obtain a category $\calV$. Finally we
define a tensor product $\otimes: \calV\times\calV\to\calV$, an
associativity constraint and a commutativity constraint which give
an AC tensor structure on $\calV$.

The most difficult step is the construction of the sets
$\Hom_\calV(X,Y)$. The sets $\Hom_\calV(X,Y)$ are defined as
$\Ho(X,Y)/\sim_{X,Y}$ where $\Ho(X,Y)$ is a set of words over an
alphabet and $\sim_{X,Y}$ is an equivalence relation on $\Ho(X,Y)$.
To ensure that $\circ$ and $\otimes$ are defined on $\calV$ and that
we have a determinant functor $f:\calT\to\calV$ we must construct
$\Ho(X,Y)$ and $\sim_{X,Y}$ for all pairs $X,Y$ simultaneously.

\textbf{Constructing the objects.} Let $\AA$ be the alphabet
\begin{equation*}
\AA = \{ \underline{A} : A \mbox{ is an object in } \calT \} \cup
\{\langle,\rangle,\boxtimes\}.
\end{equation*}
We define $\objects(\calV)$ to be the smallest subset of $\AA^*$
which satisfies the following two conditions.
\begin{itemize}
\item $\underline{A}\in\objects(\calV)$ for every object $A$ in
$\calT$.
\item If $X,Y\in\objects(\calV)$ then $\langle X\boxtimes
Y\rangle\in\objects(\calV)$.
\end{itemize}

If $X\in\objects(\calV)$ then either $X=\underline{A}$ for some
object $A$ in $\calT$, or $X=\langle Y\boxtimes Z\rangle$ for unique
$Y,Z\in\objects(\calV)$. This will allow us to argue by induction on
the length of words.

\textbf{Constructing the morphisms.} Let $\BB$ be the alphabet
\begin{equation*}
\begin{split}
\BB = & \{\iota_X : X\in\objects(\calV)\} \cup \{\langle,\rangle,\diamond,\boxtimes\} \\
& \cup \{\underline{\varphi}_{X,Y,Z}, \overline{\varphi}_{X,Y,Z} :
X,Y,Z\in\objects(\calV)\} \cup \{\underline{\psi}_{X,Y} : X,Y\in\objects(\calV)\} \\
& \cup \{ \underline{a} : a \mbox{ is an isomorphism in }\calT\}
\cup \{\underline{\Delta}, \overline{\Delta} : \Delta \mbox{ is a
distinguished triangle in }\calT\}.
\end{split}
\end{equation*}
We define the system of sets $\Ho(X,Y)\subseteq\BB^*$ (parameterized
by $X,Y\in\objects(\calV)$) to be the smallest such system which
satisfies the following nine conditions (where each condition must
hold for all $X, Y, \dots\in\objects(\calV)$).
\begin{itemize}
\item $\iota_X\in \Ho(X,X)$.
\item If $\alpha\in \Ho(X,Y)$ and $\beta\in \Ho(Y,Z)$ then
$\langle\beta\diamond\alpha\rangle\in \Ho(X,Z)$.
\item If $\alpha_1\in \Ho(X_1,Y_1)$ and $\alpha_2\in \Ho(X_2,Y_2)$ then
$\langle\alpha_1\boxtimes\alpha_2\rangle\in \Ho(\langle X_1\boxtimes
X_2\rangle,\langle Y_1\boxtimes Y_2\rangle)$.
\item $\underline{\varphi}_{X,Y,Z}\in
\Ho(\langle X\boxtimes \langle Y\boxtimes Z\rangle\rangle,
\langle\langle X\boxtimes Y\rangle\boxtimes Z\rangle)$.
\item $\overline{\varphi}_{X,Y,Z}\in \Ho(\langle\langle X\boxtimes Y\rangle\boxtimes Z\rangle,
\langle X\boxtimes\langle Y\boxtimes Z\rangle\rangle)$.
\item $\underline{\psi}_{X,Y}\in \Ho(\langle X\boxtimes Y\rangle,\langle Y\boxtimes X\rangle)$.
\item $\underline{a}\in \Ho(\underline{A},\underline{B})$ for every
isomorphism $a: A\to B$ in $\calT$.
\item $\underline{\Delta}\in \Ho(\underline{B},
\langle\underline{A}\boxtimes\underline{C}\rangle)$ for every
distinguished triangle $\Delta: A\to B\to C\to TA$ in $\calT$.
\item $\overline{\Delta}\in \Ho(\langle\underline{A}\boxtimes\underline{C}\rangle,
\underline{B})$ for every distinguished triangle $\Delta: A\to B\to
C\to TA$ in $\calT$.
\end{itemize}

If $\alpha\in \Ho(U,V)$ then either $\alpha$ is one of $\iota_X,
\underline{\varphi}_{X,Y,Z}, \overline{\varphi}_{X,Y,Z},
\underline{\psi}_{X,Y}, \underline{a}, \underline{\Delta},
\overline{\Delta}$, or $\alpha=\langle\gamma\diamond\beta\rangle$
for unique $\beta$ and $\gamma$, or
$\alpha=\langle\beta\boxtimes\gamma\rangle$ for unique $\beta$ and
$\gamma$. Hence we can argue by induction on the length of words.
Furthermore we observe that $\Ho(U,V)$ and $\Ho(U',V')$ are disjoint
unless $U=U'$ and $V=V'$.

We want to define $\calV$ to be the category whose objects are
$\objects(\calV)$ and whose morphisms from $X$ to $Y$ are
$\Ho(X,Y)/\sim_{X,Y}$, where $\sim_{X,Y}$ is an equivalence relation
on $\Ho(X,Y)$. The following definition formalizes the conditions we
must impose on $\sim_{X,Y}$.

\begin{definition}
Suppose that for each pair $X,Y\in\objects(\calV)$ we have an
equivalence relation $\sim_{X,Y}$ on $\Ho(X,Y)$. We say that such a
system $\{\sim_{X,Y}\}$ is an admissible equivalence relation if it
satisfies the following conditions (where each condition must hold
for all $X, Y, \dots\in\objects(\calV)$). To simplify the notation
we omit the subscript $X,Y$ of $\sim_{X,Y}$, and we do not
explicitly mention in which set $\Ho(X,Y)$ a relation
$\alpha\sim\beta$ holds if it is clear from the context.

\begin{itemize}

\item[(i)] The following conditions ensure that $\calV$ becomes a
category with composition $\diamond$ and identity morphisms
$\iota_X$.
\begin{itemize}
\item If $\alpha\sim\alpha'$ in $\Ho(X,Y)$ and $\beta\sim\beta'$ in $\Ho(Y,Z)$ then
$\langle\beta\diamond\alpha\rangle\sim\langle\beta'\diamond\alpha'\rangle$.
\item $\langle\iota_Y\diamond\alpha\rangle\sim\alpha$ for all $\alpha\in \Ho(X,Y)$.
\item $\langle\alpha\diamond\iota_X\rangle\sim\alpha$ for all $\alpha\in \Ho(X,Y)$.
\item $\langle\gamma\diamond\langle\beta\diamond\alpha\rangle\rangle\sim
\langle\langle\gamma\diamond\beta\rangle\diamond\alpha\rangle$ for
all $\alpha\in \Ho(X,Y)$, $\beta\in \Ho(Y,Z)$, $\gamma\in \Ho(Z,W)$.
\end{itemize}

\item[(ii)] The following conditions ensure that $\boxtimes$ becomes a
bifunctor on $\calV$.
\begin{itemize}
\item If $\alpha_1\sim\alpha'_1$ in $\Ho(X_1,Y_1)$ and $\alpha_2\sim\alpha'_2$
in $\Ho(X_2,Y_2)$ then
$\langle\alpha_1\boxtimes\alpha_2\rangle\sim\langle\alpha'_1\boxtimes\alpha'_2\rangle$.
\item $\langle\iota_X\boxtimes\iota_Y\rangle\sim\iota_{\langle X\boxtimes
Y\rangle}$.
\item If $\alpha_1\in \Ho(X_1,Y_1)$, $\beta_1\in \Ho(Y_1,Z_1)$, $\alpha_2\in
\Ho(X_2,Y_2)$, $\beta_2\in \Ho(Y_2,Z_2)$ then
$\langle\langle\beta_1\diamond\alpha_1\rangle \boxtimes
\langle\beta_2\diamond\alpha_2\rangle\rangle \sim
\langle\langle\beta_1\boxtimes\beta_2\rangle\diamond
\langle\alpha_1\boxtimes\alpha_2\rangle\rangle$.
\end{itemize}

\item[(iii)] The following conditions ensure that
$\underline{\varphi}_{X,Y,Z}$ and $\underline{\psi}_{X,Y}$ become
compatible associativity and commutativity constraints on $\calV$.
\begin{itemize}
\item ($\underline{\varphi}_{X,Y,Z}$ is natural) If $\alpha\in \Ho(X,X')$, $\beta\in
\Ho(Y,Y')$, $\gamma\in \Ho(Z,Z')$ then
$\langle\langle\langle\alpha\boxtimes\beta\rangle\boxtimes\gamma\rangle
\diamond \underline{\varphi}_{X,Y,Z}\rangle \sim
\langle\underline{\varphi}_{X',Y',Z'} \diamond
\langle\alpha\boxtimes\langle\beta\boxtimes\gamma\rangle\rangle\rangle$.
\item ($\underline{\varphi}_{X,Y,Z}$ is an isomorphism)
$\langle\underline{\varphi}_{X,Y,Z}\diamond\overline{\varphi}_{X,Y,Z}\rangle
\sim \iota_{\langle\langle X\boxtimes Y\rangle\boxtimes Z\rangle}$
and
$\langle\overline{\varphi}_{X,Y,Z}\diamond\underline{\varphi}_{X,Y,Z}\rangle
\sim \iota_{\langle X\boxtimes\langle Y\boxtimes Z\rangle\rangle}$.
\item (pentagonal axiom)
\begin{equation*}
\begin{split}
\qquad\qquad\qquad & \langle\underline{\varphi}_{\langle X\boxtimes
Y\rangle,Z,W}\diamond
\underline{\varphi}_{X,Y,\langle Z\boxtimes W\rangle}\rangle \\
& \qquad\qquad\qquad \sim
\langle\langle\underline{\varphi}_{X,Y,Z}\boxtimes
\iota_W\rangle\diamond \langle\underline{\varphi}_{X,Y\boxtimes
Z,W}\diamond \langle\iota_X\boxtimes
\underline{\varphi}_{Y,Z,W}\rangle\rangle\rangle.
\end{split}
\end{equation*}
\item ($\underline{\psi}_{X,Y}$ is natural) If $\alpha\in \Ho(X,X')$ and
$\beta\in \Ho(Y,Y')$ then
$\langle\langle\beta\boxtimes\alpha\rangle\diamond\underline{\psi}_{X,Y}\rangle
\sim \langle\underline{\psi}_{X',Y'}\diamond
\langle\alpha\boxtimes\beta\rangle\rangle$.
\item $\langle\underline{\psi}_{Y,X}\diamond\underline{\psi}_{X,Y}\rangle\sim
\iota_{\langle X\boxtimes Y\rangle}$.
\item (hexagonal axiom)
\begin{equation*}
\begin{split}
\qquad\qquad & \langle\underline{\varphi}_{Z,X,Y}\diamond
\langle\underline{\psi}_{\langle X\boxtimes Y\rangle,Z}\diamond
\underline{\varphi}_{X,Y,Z}\rangle\rangle \\
& \qquad\qquad\qquad\qquad \sim
\langle\langle\underline{\psi}_{X,Z}\boxtimes \iota_Y\rangle\diamond
\langle\underline{\varphi}_{X,Z,Y}\diamond
\langle\iota_X\boxtimes\underline{\psi}_{Y,Z}\rangle\rangle\rangle.
\end{split}
\end{equation*}
\end{itemize}

\item[(iv)] The following conditions ensure that we obtain a functor $f_1:
\calT_{\iso}\to\calV$ and isomorphisms $f_2(\Delta)$ in $\calV$.
\begin{itemize}
\item $\underline{\id_A}\sim\iota_{\underline{A}}$ for every object $A$ of $\calT$.
\item $\underline{b\circ a}\sim \langle\underline{b}\diamond\underline{a}\rangle$
for all isomorphisms $a: A\to B$, $b: B\to C$ in $\calT$.
\item ($\underline{\Delta}$ is an isomorphism)
$\langle\underline{\Delta}\diamond\overline{\Delta}\rangle \sim
\iota_{\langle\underline{A}\boxtimes\underline{C}\rangle}$ and
$\langle\overline{\Delta}\diamond\underline{\Delta}\rangle \sim
\iota_{\underline{B}}$ for every distinguished triangle $\Delta:
A\to B\to C\to TA$ in $\calT$.
\end{itemize}

\item[(v)] The following conditions ensure that the axioms for a
determinant functor hold.
\begin{itemize}
\item (naturality) $\langle\underline{\Delta'}\diamond\underline{b}\rangle \sim
\langle\langle\underline{a}\boxtimes\underline{c}\rangle\diamond\underline{\Delta}\rangle$
for every isomorphism of distinguished triangles as in
(\ref{diagram_iso_dist}).
\item (associativity)
$\langle\underline{\varphi}_{\underline{A},\underline{C'},\underline{A'}}
\diamond
\langle\langle\iota_{\underline{A}}\boxtimes\underline{\Delta_{\mathrm{v2}}}\rangle
\diamond \underline{\Delta_{\mathrm{h2}}}\rangle\rangle \sim
\langle\langle\underline{\Delta_{\mathrm{h1}}}\boxtimes\iota_{\underline{A'}}\rangle
\diamond \underline{\Delta_{\mathrm{v1}}}\rangle$ for every
octahedral diagram as in (\ref{eqn_octahedral_diagram}) (where
$\Delta_{\mathrm{h1}}$ and $\Delta_{\mathrm{h2}}$ (resp.\
$\Delta_{\mathrm{v1}}$ and $\Delta_{\mathrm{v2}}$) are the first and
second horizontal (resp.\ vertical) distinguished triangles in
diagram (\ref{eqn_octahedral_diagram})).
\item (commutativity)
$\langle\underline{\psi}_{\underline{A},\underline{B}} \diamond
\underline{\Delta_1}\rangle \sim \underline{\Delta_2}$ for every
pair of distinguished triangles as in (\ref{eqn_split_dist}).
\end{itemize}

\end{itemize}
\end{definition}

The intersection of admissible equivalence relations is again an
admissible equivalence relation. Therefore there exists a unique
minimal admissible equivalence relation which from now on we will
denote by $\{\sim_{X,Y}\}$. For every pair of objects $X,Y$ in
$\calV$ we define $\Hom_{\calV}(X,Y)=\Ho(X,Y)/\sim_{X,Y}$. We denote
the equivalence class of $\alpha\in \Ho(X,Y)$ in
$\Hom_\calV(X,Y)=\Ho(X,Y)/\sim_{X,Y}$ by $[\alpha]$.

\textbf{Constructing the Picard category.} We define $\calV$ to be
the following category: the set of objects of $\calV$ is
$\objects(\calV)$, the set of morphisms from $X\in\objects(\calV)$
to $Y\in\objects(\calV)$ is $\Hom_{\calV}(X,Y)$, and the composition
of morphisms $[\alpha]\in\Hom_{\calV}(X,Y)$ and
$[\beta]\in\Hom_{\calV}(Y,Z)$ is $[\beta]\circ [\alpha] =
[\langle\beta\diamond\alpha\rangle]\in\Hom_\calV(X,Z)$. The
definition of the equivalence relation ensures that this composition
is well-defined and that we obtain a category. The identity morphism
$\id_X$ of an object $X$ is given by $\id_X=[\iota_X]$.

Next we define a bifunctor $\otimes: \calV\times\calV\to\calV$ by
$X\otimes Y=\langle X\boxtimes Y\rangle$ for objects $X,Y$ and
$[\alpha]\otimes [\beta]=[\langle\alpha\boxtimes\beta\rangle]$ for
morphisms $[\alpha],[\beta]$. Then
$\varphi_{X,Y,Z}=[\underline{\varphi}_{X,Y,Z}] \in
\Hom_\calV(X\otimes(Y\otimes Z),(X\otimes Y)\otimes Z)$ and
$\psi_{X,Y}=[\underline{\psi}_{X,Y}]\in\Hom_\calV(X\otimes
Y,Y\otimes X)$ give compatible associativity and commutativity
constraints, so that $\calV$ becomes an AC tensor category.

\begin{lemma}
\label{lemma_V_is_Picard} The AC tensor category $\calV$ is a Picard
category.
\end{lemma}

\begin{proof}
Clearly $\calV$ is non-empty. Also, since each of $[\iota_X],
[\underline{\varphi}_{X,Y,Z}], [\overline{\varphi}_{X,Y,Z}]$,
$[\underline{\psi}_{X,Y}]$, $[\underline{a}]$,
$[\underline{\Delta}]$, $[\overline{\Delta}]$ is an isomorphism, it
follows by induction that every morphism in $\calV$ is an
isomorphism.

Next we show that for every object $X$ in $\calV$ there exists an
object $A$ in $\calT$ and an isomorphism $[\alpha]:
X\to\underline{A}$ in $\calV$ such that $\alpha\in
\Ho(X,\underline{A})$ is a word which does not contain any
$\underline{\psi}$. If $X=\underline{A}$ we can take
$\alpha=\iota_{\underline{A}}$. Now suppose that $X=\langle
Y\boxtimes Z\rangle$ and that there exist isomorphisms $[\beta]:
Y\to\underline{B}$ and $[\gamma]: Z\to\underline{C}$ where $B$ and
$C$ are objects in $\calT$ and $\beta\in \Ho(Y,\underline{B})$ and
$\gamma\in \Ho(Z,\underline{C})$ do not contain any
$\underline{\psi}$. Let $\Delta$ be a distinguished triangle of the
form $B\to A\to C\to TB$. Then we can take
$\alpha=\langle\overline{\Delta}\diamond
\langle\beta\boxtimes\gamma\rangle\rangle\in \Ho(X,\underline{A})$.

Now fix a zero object $0$ in $\calT$. We want to show that there
exists a functorial isomorphism $\lambda_X: \underline{0}\otimes
X\to X$ in $\calV$. For $X=\underline{A}$ with $A$ an object in
$\calT$ we let $\lambda_X=[\overline{\Delta}]$ where $\Delta: 0\to
A\xrightarrow{\id} A\to T0$. For $X, Y\in\objects(\calV)$ we define
$\lambda_{\langle X\boxtimes Y\rangle}$ as the composite morphism
\begin{equation*}
\underline{0}\otimes\langle X\boxtimes Y\rangle
=\underline{0}\otimes(X\otimes Y)
\xrightarrow{\varphi_{\underline{0},X,Y}} (\underline{0}\otimes
X)\otimes Y\xrightarrow{\lambda_X\otimes\id_Y} X\otimes Y = \langle
X\boxtimes Y\rangle.
\end{equation*}

To show that $\lambda_X$ is functorial we must verify that
\begin{equation}
\label{eqn_lambda_functorial}
[\alpha]\circ\lambda_X=\lambda_Y\circ(\id_{\underline{0}}\otimes[\alpha])
\end{equation}
in $\calV$ for all $\alpha\in \Ho(X,Y)$. For this we use induction
on the length of $\alpha$. It is easy to check
(\ref{eqn_lambda_functorial}) for $\alpha=\iota_X$ and
$\alpha=\underline{a}$. If $\alpha=\underline{\Delta}$ for a
distinguished triangle $\Delta: A\to B\to C\to TA$ then
(\ref{eqn_lambda_functorial}) follows from the associativity axiom
applied to the octahedral diagram
\begin{equation*}
\xymatrix{ 0 \ar[r] \ar@{=}[d] & A \ar[r]^{\id} \ar[d] & A \ar[r]
\ar[d] & T0 \ar@{=}[d] \\
0 \ar[r] & B \ar[r]^{\id} \ar[d] & B \ar[r] \ar[d] & T0 \\
& C \ar@{=}[r] \ar[d] & C \ar[d] & \\
& TA \ar[r]^{\id} & TA. & }
\end{equation*}
The statement also follows for $\alpha=\overline{\Delta}$ because
$[\overline{\Delta}]=[\underline{\Delta}]^{-1}$ in $\calV$. If
$\alpha=\underline{\varphi}_{X,Y,Z}$ then we can deduce
(\ref{eqn_lambda_functorial}) from the pentagonal axiom and
naturality of $\underline{\varphi}_{X,Y,Z}$, and if
$\alpha=\overline{\varphi}_{X,Y,Z}$ then we use that
$[\overline{\varphi}_{X,Y,Z}]=[\underline{\varphi}_{X,Y,Z}]^{-1}$ in
$\calV$. Next if (\ref{eqn_lambda_functorial}) holds for $\alpha$
and $\beta$ then it holds for $\langle\beta\diamond\alpha\rangle$
(if defined) and for $\langle\alpha\boxtimes\beta\rangle$. Hence it
follows by induction that (\ref{eqn_lambda_functorial}) is true for
all $\alpha$ that do not contain any $\underline{\psi}$. Now
consider $\alpha=\underline{\psi}_{X,Y}$. We have shown above that
there exist objects $A, B$ in $\calT$ and isomorphisms
$X\cong\underline{A}$, $Y\cong\underline{B}$ in $\calV$ such that
these isomorphisms are represented by words in
$\Ho(X,\underline{A})$ and $\Ho(Y,\underline{B})$ which do not
contain any $\underline{\psi}$. Hence for these isomorphisms we know
(\ref{eqn_lambda_functorial}) which by naturality of
$\underline{\psi}$ reduces to the case
$\alpha=\underline{\psi}_{\underline{A},\underline{B}}$. In this
case (\ref{eqn_lambda_functorial}) follows from the commutativity
axiom and the case of $\underline{\Delta}$ already shown above.
Finally, using induction again we can deduce the validity of
(\ref{eqn_lambda_functorial}) for all $\alpha\in \Ho(X,Y)$.

Now let $W\in\objects(\calV)$. We must show that the functor
$X\mapsto W\otimes X$ is an autoequivalence of $\calV$. Since
$W\cong \underline{A}$ for some object $A$ in $\calT$, it suffices
to show this for the functor $F: X\mapsto \underline{A}\otimes X$.
Let $G$ be the functor $Y\mapsto \underline{TA}\otimes Y$. Then
$G\circ F$ is isomorphic to the identity functor by
\begin{equation*}
\underline{TA}\otimes (\underline{A}\otimes X)\cong
(\underline{TA}\otimes \underline{A})\otimes X\cong
(\underline{A}\otimes \underline{TA})\otimes
X\xrightarrow{[\overline{\Delta}]\otimes\id} \underline{0}\otimes
X\xrightarrow{\lambda_X} X
\end{equation*}
where $\Delta: A\to 0\to TA\xrightarrow{\id} TA$, and $F\circ G$ is
isomorphic to the identity functor by
\begin{equation*}
\underline{A}\otimes (\underline{TA}\otimes X)\cong
(\underline{A}\otimes \underline{TA})\otimes
X\xrightarrow{[\overline{\Delta}]\otimes\id} \underline{0}\otimes
X\xrightarrow{\lambda_X} X.
\end{equation*}
This completes the proof of Lemma \ref{lemma_V_is_Picard}.
\end{proof}

\textbf{Constructing the determinant functor.} For an object $A$ in
$\calT$ we define $f_1(A)=\underline{A}\in\objects(\calV)$, and for
an isomorphism $a: A\to B$ in $\calT$ we define
$f_1(a)=[\underline{a}]\in\Hom_{\calV}(f_1(A),f_1(B))$. For a
distinguished triangle $\Delta: A\to B\to C\to TA$ we define
$f_2(\Delta)=[\underline{\Delta}]\in\Hom_{\calV}(f_1(B),f_1(A)\otimes
f_1(C))$.

\begin{lemma}
$f=(f_1,f_2)$ is a determinant functor on $\calT$ with values in
$\calV$.
\end{lemma}

\begin{proof}
This is immediate from the definition of the equivalence relation.
\end{proof}

\textbf{Verification of the universal property.} We must show that
for every Picard category $\calP$ the functor
$\Hom^{\otimes}(\calV,\calP)\to\det(\calT,\calP)$, $M\mapsto M\circ
f$, is an equivalence of categories.

\begin{lemma}
\label{lemma_existence_2} The functor
$\Hom^{\otimes}(\calV,\calP)\to\det(\calT,\calP)$ is surjective on
objects.
\end{lemma}

\begin{proof}
Let $g=(g_1,g_2):\calT\to\calP$ be a determinant functor. We want to
show that there exists a monoidal functor $M:\calV\to\calP$ such
that $g=M\circ f$.

On objects of $\calV$ we define $M$ by $M(\underline{A})=g_1(A)$ and
$M(\langle X\boxtimes Y\rangle) = M(X)\otimes M(Y)$. Next we define
a morphism $M(\alpha): M(X)\to M(Y)$ for every $\alpha\in \Ho(X,Y)$.
Firstly $M(\iota_X)=\id_{M(X)}$. If $a: A\to B$ is an isomorphism in
$\calT$ then $M(\underline{a})=g_1(a)$. If $\Delta: A\to B\to C\to
TA$ is a distinguished triangle in $\calT$ then
$M(\underline{\Delta})=g_2(\Delta)$ and
$M(\overline{\Delta})=g_2(\Delta)^{-1}$. Also
$M(\underline{\varphi}_{X,Y,Z})=\varphi_{M(X),M(Y),M(Z)}$,
$M(\overline{\varphi}_{X,Y,Z})=\varphi^{-1}_{M(X),M(Y),M(Z)}$ and
$M(\underline{\psi}_{X,Y})=\psi_{M(X),M(Y)}$. Finally
$M(\langle\beta\diamond\alpha\rangle)=M(\beta)\circ M(\alpha)$ and
$M(\langle\alpha_1\boxtimes\alpha_2\rangle) = M(\alpha_1)\otimes
M(\alpha_2)$.

For a morphism $[\alpha]\in\Hom_\calV(X,Y)=\Ho(X,Y)/\sim_{X,Y}$ we
now define $M([\alpha])=M(\alpha)\in\Hom_\calP(M(X),M(Y))$. To show
that this is well-defined we must verify that the equivalence
relation defined by $\alpha\sim\alpha'$ if $M(\alpha)=M(\alpha')$ is
admissible. This is straightforward to check using that $\calP$ is a
Picard category and $g=(g_1,g_2)$ is a determinant functor.

Finally we define a natural isomorphism $c_{X,Y}: M(X)\otimes
M(Y)\to M(X\otimes Y)$ by $c_{X,Y}=\id_{M(X)\otimes M(Y)}$. Then it
is easy to check that $M=(M,c)$ is a monoidal functor
$\calV\to\calP$ and that $g=M\circ f$.
\end{proof}

\begin{lemma}
\label{lemma_existence_3} The functor
$\Hom^{\otimes}(\calV,\calP)\to\det(\calT,\calP)$ is fully faithful.
\end{lemma}

\begin{proof}
Let $M=(M,c), M'=(M',c'):\calV\to\calP$ be two monoidal functors and
let $\lambda: M\circ f \to M'\circ f$ be an isomorphism in
$\det(\calT,\calP)$. We must show that there exists a unique
isomorphism $\kappa: M\to M'$ in $\Hom^{\otimes}(\calV,\calP)$ such
that $\lambda=\kappa\circ f$.

If there is an isomorphism $\kappa: M\to M'$ such that
$\lambda=\kappa\circ f$ then necessarily
$\kappa_{\underline{A}}=\lambda_A$ for every object $A$ in $\calT$
and $\kappa_{\langle X\boxtimes Y\rangle} =
(c'_{X,Y})\circ(\kappa_X\otimes\kappa_Y)\circ c_{X,Y}^{-1}$ for all
$X,Y\in\objects(\calV)$. This shows that $\kappa$ is uniquely
determined by $\lambda$.

To show the existence of $\kappa$ we define a morphism $\kappa_X:
M(X)\to M'(X)$ for every object $X$ in $\calV$ as follows:
$\kappa_{\underline{A}}=\lambda_A$ for every object $A$ of $\calT$,
$\kappa_{\langle X\boxtimes Y\rangle}=(c'_{X,Y})\circ
(\kappa_X\otimes\kappa_Y)\circ c_{X,Y}^{-1}$. Using induction on the
length of words which represent morphisms in $\calV$, one can easily
check that $\kappa$ is an isomorphism of monoidal functors $M\to
M'$, and by definition of $\kappa$ it is obvious that
$\lambda=\kappa\circ f$.
\end{proof}

It follows from Lemmas \ref{lemma_existence_2} and
\ref{lemma_existence_3} that for every Picard category $\calP$ the
functor $\Hom^{\otimes}(\calV,\calP)\to\det(\calT,\calP)$ is an
equivalence of categories. This completes the proof of Theorem
\ref{theorem_existence}.
\end{proof}

\subsection{$K_0$ and $K_1$ of a triangulated category}

\begin{definition}
Let $\calT$ be a small triangulated category. For $i=0, 1$ we define
the group $K_i(\calT)$ by $K_i(\calT)=\pi_i(\calV)$ where $f:
\calT\to\calV$ is a universal determinant functor on $\calT$.
\end{definition}

The following argument shows that the groups $K_i(\calT)$ are
well-defined and functorial in $\calT$. Let $F: \calS\to\calT$ be an
exact functor of small triangulated categories and $f_{\calS}:
\calS\to\calV(\calS)$, $f_{\calT}: \calT\to\calV(\calT)$ universal
determinant functors for $\calS$ and $\calT$ respectively. The
equivalence of categories
$\Hom^\otimes(\calV(\calS),\calV(\calT))\to\det(\calS,\calV(\calT))$
shows that there exists a monoidal functor
$M:\calV(\calS)\to\calV(\calT)$, unique up to isomorphism, such that
the determinant functors $M\circ f_{\calS}$ and $f_{\calT}\circ F$
in $\det(\calS,\calV(\calT))$ are isomorphic. The monoidal functor
$M$ induces homomorphisms $\pi_i(M):
\pi_i(\calV(\calS))\to\pi_i(\calV(\calT))$ for $i=0,1$, and in fact
these homomorphisms only depend on $F$ because isomorphic monoidal
functors induce the same maps on $\pi_i$. Applying this to the
identity functor $F=\id: \calT\to\calT$ implies that up to canonical
isomorphism $K_i(\calT)$ is independent of the choice of universal
determinant functor. For an arbitrary exact functor $F:
\calS\to\calT$ one obtains an induced homomorphism $K_i(F)=\pi_i(M):
K_i(\calS)\to K_i(\calT)$.

\begin{lemma}
Let $F, F': \calS\to\calT$ be exact functors of small triangulated
categories. If $F$ and $F'$ are isomorphic then $K_i(F)=K_i(F'):
K_i(\calS)\to K_i(\calT)$ for $i=0,1$.
\end{lemma}

\begin{proof}
Let $f_\calS:\calS\to\calV(\calS)$ and
$f_\calT:\calT\to\calV(\calT)$ be universal determinant functors.
Since $F$ and $F'$ are isomorphic exact functors it follows that
$f_\calT\circ F$ and $f_\calT\circ F'$ are isomorphic determinant
functors. Hence there exists a monoidal functor
$M:\calV(S)\to\calV(\calT)$ such that $M\circ f_\calS\cong
f_\calT\circ F$ and $M\circ f_\calS\cong f_\calT\circ F'$. By
definition $K_i(F)=\pi_i(M)$ and $K_i(F')=\pi_i(M)$.
\end{proof}

This lemma immediately implies

\begin{corollary}
\label{corollary_K_equivalence} If $F:\calS\to\calT$ is an
equivalence of small triangulated categories then $K_i(F):
K_i(\calS)\to K_i(\calT)$ is an isomorphism.
\end{corollary}

\begin{remark}
Let $\calT$ be a triangulated category such that the isomorphism
classes of objects in $\calT$ form a set. Then there exists a small
full subcategory $\calS$ of $\calT$ satisfying $T\calS=\calS$ and
for which the inclusion functor $\calS\to\calT$ is an equivalence of
categories. The triangulated category $\calT$ induces the structure
of a triangulated category on $\calS$ and we can define $K_i(\calT)$
for $i=0,1$ by $K_i(\calT)=K_i(\calS)$. Corollary
\ref{corollary_K_equivalence} shows that up to canonical isomorphism
the definition of $K_i(\calT)$ does not depend on the choice of
$\calS$.
\end{remark}

Next we show that the group $K_0(\calT)$ is the Grothendieck group
of the triangulated category $\calT$. Recall that the Grothendieck
group is defined to be the quotient $F(\calT)/R(\calT)$ where
$F(\calT)$ is the free abelian group on isomorphism classes of
objects of $\calT$, and $R(\calT)$ is the subgroup of $F(\calT)$
which is generated by all elements $(B)-(A)-(C)$ whenever there is a
distinguished triangle $A\to B\to C\to TA$ (here $(A)$ denotes the
element in $F(\calT)$ corresponding to the isomorphism class of the
object $A$).

\begin{proposition}
The group $K_0(\calT)$ is canonically isomorphic to the Grothendieck
group of $\calT$.
\end{proposition}

\begin{proof}
Let $f=(f_1,f_2): \calT\to\calV$ be a universal determinant functor
for $\calT$. We define a homomorphism $F(\calT)\to K_0(\calT)$ by
sending $(A)$ to the isomorphism class of $f_1(A)$ in
$\pi_0(\calV)=K_0(\calT)$. One easily verifies that this
homomorphism factors through a homomorphism $u:F(\calT)/R(\calT)\to
K_0(\calT)$.

Let $\calP$ be the Picard category whose set of objects is the set
$F(\calT)/R(\calT)$, whose morphisms are only the identity morphism
of each object and whose tensor product is given by the product in
the group $F(\calT)/R(\calT)$. There is a unique determinant functor
$g: \calT\to\calP$ which sends an object $A$ of $\calT$ to the image
of $(A)$ in $F(\calT)/R(\calT)$. By the universality of $f:
\calT\to\calV$ there exists a monoidal functor $M: \calV\to\calP$
(unique up to isomorphism) such that $M\circ f\cong g$, and by
applying $\pi_0$ we obtain a homomorphism $v:
K_0(\calT)=\pi_0(\calV) \xrightarrow{\pi_0(M)}
\pi_0(\calP)=F(\calT)/R(\calT)$.

It is not difficult to verify that $vu=\id$ and that $uv(X)=X$ if
$X\in\pi_0(\calV)=K_0(\calT)$ is the isomorphism class of $f_1(A)$
for some object $A$ in $\calT$. But the proof of Theorem
\ref{theorem_existence} shows that every object in $\calV$ is
isomorphic to $f_1(A)$ for some $A$, hence $uv=\id$ follows.
\end{proof}

We conclude this section with a description of the automorphisms of
a determinant functor on $\calT$ in terms of homomorphisms on
$K_0(\calT)$.

\begin{lemma}
Let $\calT$ be a small triangulated category, $\calP$ a Picard
category and $g: \calT\to\calP$ a determinant functor. Then there
exists a canonical isomorphism
$\Aut_{\det(\calT,\calP)}(g)\cong\Hom(K_0(\calT),\pi_1(\calP))$
(where $\Hom$ is homomorphisms of abelian groups).
\end{lemma}

\begin{proof}
Let $f:\calT\to\calV$ be a universal determinant functor and
$M:\calV\to\calP$ a monoidal functor for which $M\circ f\cong g$ in
$\det(\calT,\calP)$. Then $\Aut_{\det(\calT,\calP)}(g)\cong
\Aut_{\det(\calT,\calP)}(M\circ f)\cong
\Aut_{\Hom^\otimes(\calV,\calP)}(M)\cong
\Hom(\pi_0(\calV),\pi_1(\calP))= \Hom(K_0(\calT),\pi_1(\calP))$,
where the second isomorphism follows from the universality of $f$
and the third isomorphism is Lemma
\ref{lemma_monoidal_automorphism}.
\end{proof}

\section{Determinant functors on triangulated categories with t-structure}
\label{section_t_structure}

In this section we study the relation between the determinant
functors on a triangulated category with a t-structure and the
determinant functors on its heart.

\subsection{Statement of result}

Let $\calT$ be a triangulated category with t-structure and let
$\calE$ be its heart. For the precise definition and a detailed
discussion of these notions we refer the reader to \cite[\S
1.3]{BBD}. Here we just recall that a t-structure on a triangulated
category $\calT$ is a pair of full subcategories $(\calT^{\leq
0},\calT^{\geq 0})$ of $\calT$ satisfying various properties, and
that the heart $\calE$ of the t-structure is the intersection
$\calE=\calT^{\leq 0}\cap\calT^{\geq 0}$. The heart is a full
abelian subcategory of $\calT$. If $\Delta: 0\to A\xrightarrow{a}
B\xrightarrow{b} C\to 0$ is a short exact sequence in $\calE$ then
there exists a unique morphism $c: C\to TA$ in $\calT$ such that
$A\xrightarrow{a} B\xrightarrow{b} C\xrightarrow{c} TA$ is a
distinguished triangle in $\calT$. We will denote this distinguished
triangle by $\widehat{\Delta}$.

Let $\calP$ be a Picard category and $f=(f_1,f_2)$ a determinant
functor on $\calT$ with values in $\calP$. We define a functor $g_1:
\calE_\iso\to\calP$ to be the restriction of $f_1$ to $\calE_\iso$,
and for every short exact sequence $\Delta: 0\to A\to B\to C\to 0$
in $\calE$ we define an isomorphism $g_2(\Delta): g_1(B)\to
g_1(A)\otimes g_1(C)$ by $g_2(\Delta)=f_2(\widehat{\Delta})$. Let
$f|_{\calE}$ denote the pair $(g_1,g_2)$.

\begin{lemma}
\label{lemma_restriction} If $f$ is a determinant functor on the
triangulated category $\calT$ with values in $\calP$, then
$f|_{\calE}$ is a determinant functor on the exact category $\calE$
with values in $\calP$. Moreover $f\mapsto f|_{\calE}$ naturally
extends to a functor $\det(\calT,\calP)\to\det(\calE,\calP)$.
\end{lemma}

\begin{proof}
If $\Delta\to\Delta'$ is an isomorphism of short exact sequences in
$\calE$ then we obtain an induced isomorphism of distinguished
triangles $\widehat{\Delta}\to \widehat{\Delta'}$ in $\calT$. Hence
the naturality axiom for $f$ implies the naturality axiom for
$f|_{\calE}$. Given a diagram of short exact sequences as in
(\ref{diagram_assoc_exact}) then the associated diagram of
distinguished triangles is an octahedral diagram, and therefore the
associativity axiom for $f|_{\calE}$ follows from the associativity
axiom for $f$. If $\Delta$ is the short exact sequence $A\to A\oplus
B\to B$ then $\widehat{\Delta}$ is the distinguished triangle $A\to
A\oplus B\to B\xrightarrow{0} TA$. Therefore the commutativity axiom
for $f|_{\calE}$ follows from the commutativity axiom for $f$.

We obtain a functor $\det(\calT,\calP)\to\det(\calE,\calP)$ because
a morphism $\lambda: f\to f'$ in $\det(\calT,\calP)$ restricts to a
morphism $\lambda|_{\calE}: f|_{\calE}\to f'|_{\calE}$ in
$\det(\calE,\calP)$.
\end{proof}

Associated to the t-structure there is a canonical cohomological
functor $H: \calT\to\calE$. The t-structure is called bounded if for
every object $A$ of $\calT$ the cohomology objects $H(T^iA)$ are
zero for all but finitely many $i\in\ZZ$. For any $n\in\ZZ$ we write
$\calT^{\leq n}=T^{-n}\calT^{\leq 0}$ and $\calT^{\geq
n}=T^{-n}\calT^{\geq 0}$. The t-structure is called non-degenerate
if the intersections $\bigcap_{n\in\ZZ}\calT^{\leq n}$ and
$\bigcap_{n\in\ZZ}\calT^{\geq n}$ consist only of zero objects. We
can now state the main result of \S \ref{section_t_structure}.

\begin{theorem}
\label{theorem_t_structure} Let $\calT$ be a triangulated category
with a non-degenerate and bounded t-structure, and let $\calE$ be
its heart. Then for every Picard category $\calP$ the restriction
functor $\det(\calT,\calP)\to\det(\calE,\calP)$, $f\mapsto
f|_{\calE}$, is an equivalence of categories.
\end{theorem}

\begin{corollary}
\label{corollary_heart} Let $\calT$ be a small triangulated category
with a non-degenerate and bounded t-structure, and let $\calE$ be
its heart. Then there exist canonical isomorphisms $K_i(\calT)\cong
K_i(\calE)$ for $i=0, 1$.
\end{corollary}

\begin{proof}
Fix universal determinant functors $f_\calT: \calT\to\calV(\calT)$
and $f_\calE: \calE\to\calV(\calE)$. Then by definition of a
universal determinant functor there exists a monoidal functor $M:
\calV(\calE)\to\calV(\calT)$ (unique up to isomorphism) such that
the determinant functors $(f_\calT)|_\calE$ and $M\circ f_\calE$ are
isomorphic. For every Picard category $\calP$ the functor
$\Hom^\otimes(\calV(\calT),\calP)\to\Hom^\otimes(\calV(\calE),\calP)$,
$N\mapsto N\circ M$, is an equivalence of categories, because the
diagram
\begin{equation*}
\xymatrix{ \Hom^\otimes(\calV(\calT),\calP) \ar[r] \ar[d] &
\Hom^\otimes(\calV(\calE),\calP) \ar[d] \\
\det(\calT,\calP) \ar[r] & \det(\calE,\calP) }
\end{equation*}
commutes up to natural isomorphism and the two vertical and the
bottom horizontal functors are equivalences (by definition of a
universal determinant functor and by Theorem
\ref{theorem_t_structure} respectively). It follows that $M$ is an
equivalence of Picard categories and hence induces isomorphisms
$K_i(\calE)=\pi_i(\calV(\calE)) \xrightarrow{\pi_i(M)}
\pi_i(\calV(\calT))=K_i(\calT)$ for $i=0,1$.
\end{proof}

As an important special case we note that if $\calE$ is an abelian
category and $\calT=\calD^\rmb(\calE)$ its bounded derived category
with the canonical t-structure then the assumptions of Corollary
\ref{corollary_heart} are satisfied and hence
$K_i(\calD^\rmb(\calE))\cong K_i(\calE)$ for $i=0,1$.

To prove Theorem \ref{theorem_t_structure} we factorize the
restriction functor $\det(\calT,\calP)\to\det(\calE,\calP)$ as the
composite of two functors
\begin{equation}
\label{equation_functor_factorization} \det(\calT,\calP)
\xrightarrow{V} \det\big((\Gr^\rmb(\calE),T),\calP\big)
\xrightarrow{W} \det(\calE,\calP)
\end{equation}
and show that each of the functors $V$ and $W$ is an equivalence of
categories.

The pair $(\Gr^\rmb(\calE),T)$ consists of the category
$\Gr^\rmb(\calE)$ of bounded $\ZZ$-graded objects in $\calE$ and the
canonical translation functor $T$ on $\Gr^\rmb(\calE)$. We will
define determinant functors on exact categories with translation in
\S \ref{subsection_exact_translation}, where we will also show that
for any exact category $\calE$ the natural functor $W:
\det\big((\Gr^\rmb(\calE),T),\calP\big) \to \det(\calE,\calP)$ is an
equivalence of categories. This deals with the functor $W$ in
(\ref{equation_functor_factorization}).

Next in \S \ref{subsection_graded_restriction} we define the functor
$V: \det(\calT,\calP) \to \det\big((\Gr^\rmb(\calE),T),\calP\big)$
for any triangulated category $\calT$ with t-structure and heart
$\calE$, and we show that the restriction functor
$\det(\calT,\calP)\to\det(\calE,\calP)$, $f\mapsto f|_{\calE}$, is
the composite $W\circ V$.

In \S \ref{subsection_cohomological_functor} we show that a bounded
cohomological functor $H$ from a triangulated category $\calT$ to an
abelian category $\calE$ induces a functor
\begin{equation*}
H^*: \det\big((\Gr^\rmb(\calE),T),\calP\big) \to \det(\calT,\calP).
\end{equation*}

Finally in \S \ref{subsection_cohomology} we complete the proof of
Theorem \ref{theorem_t_structure} by showing that if the t-structure
on $\calT$ is non-degenerate and bounded then $V$ is an equivalence
of categories. More precisely, we will show that under these
assumptions the functor $V$ is a quasi-inverse of the functor $H^*$
where $H$ is the canonical cohomological functor from $\calT$ to its
heart $\calE$.

\subsection{Exact categories with translation}
\label{subsection_exact_translation}

If $\calE$ is an exact category and $T$ an exact automorphism of
$\calE$ then we will refer to $T$ as a translation functor on
$\calE$ and to the pair $(\calE,T)$ as an exact category with
translation.

\begin{definition}
Let $(\calE,T)$ be an exact category with translation and $\calP$ a
Picard category. A determinant functor on $(\calE,T)$ with values in
$\calP$ is a triple $(g_1,g_2,\mu)$, where $(g_1,g_2)$ is a
determinant functor on the exact category $\calE$ with values in
$\calP$ and, for every object $A$ in $\calE$, $\mu_A$ is a unit
structure on $g_1(A)\otimes g_1(TA)$, i.e.\ an isomorphism
\begin{equation*}
\mu_A: g_1(A)\otimes g_1(TA)\to\big(g_1(A)\otimes
g_1(TA)\big)\otimes\big(g_1(A)\otimes g_1(TA)\big),
\end{equation*}
such that the following two conditions hold.
\begin{itemize}
\item[(i)] If $a: A\to B$ is an isomorphism in $\calE$
then $g_1(a)\otimes g_1(Ta)$ is an isomorphism of units
$(g_1(A)\otimes g_1(TA),\mu_A)\to (g_1(B)\otimes g_1(TB),\mu_B)$.
\item[(ii)] If $\Delta: 0\to A\xrightarrow{a} B\xrightarrow{b} C\to 0$ is a
short exact sequence in $\calE$ and $\Delta'$ denotes the short
exact sequence $0\to TA\xrightarrow{Ta} TB\xrightarrow{Tb} TC\to 0$,
then the composite map
\begin{equation*}
\begin{split}
\quad g_1(B)\otimes g_1(TB) & \xrightarrow{g_2(\Delta)\otimes
g_2(\Delta')}(g_1(A)\otimes g_1(C))\otimes (g_1(TA)\otimes g_1(TC)) \\
& \xrightarrow{\,\cong\,} (g_1(A)\otimes g_1(TA))\otimes
(g_1(C)\otimes g_1(TC))
\end{split}
\end{equation*}
is an isomorphism of units $(g_1(B)\otimes g_1(TB),\mu_B)\to
(g_1(A)\otimes g_1(TA),\mu_A)\otimes (g_1(C)\otimes g_1(TC),\mu_C)$.
\end{itemize}

A morphism $\lambda: (g_1,g_2,\mu)\to(g'_1,g'_2,\mu')$ of
determinant functors on $(\calE,T)$ with values in $\calP$ is a
morphism $\lambda: (g_1,g_2)\to (g'_1,g'_2)$ of the determinant
functors on the underlying exact category which is compatible with
the unit structures $\mu$ and $\mu'$, i.e.\ for every object $A$ in
$\calE$ the map $\lambda_A\otimes\lambda_{TA}$ is an isomorphism of
units $(g_1(A)\otimes g_1(TA),\mu_A)\to (g'_1(A)\otimes
g'_1(TA),\mu'_A)$.
\end{definition}

We denote the category of determinant functors on $(\calE,T)$ with
values in $\calP$ by $\det\big((\calE,T),\calP\big)$.

If $\calE$ is an exact category, then we write
$\Gr^\mathrm{b}(\calE)$ for the exact category whose objects are the
bounded $\ZZ$-graded objects in $\calE$, i.e.\ tuples
$(A_i)_{i\in\ZZ}$ of objects in $\calE$ where $A_i$ is a zero object
for all but finitely many $i$, and morphisms $(A_i)_{i\in\ZZ}\to
(B_i)_{i\in\ZZ}$ are tuples $(a_i)_{i\in\ZZ}$ where each $a_i:
A_i\to B_i$ is a morphism in $\calE$. The category
$\Gr^\mathrm{b}(\calE)$ has a canonical translation functor $T$
given by $T((A_i)_{i\in\ZZ})_j=A_{j+1}$ for objects and
$T((a_i)_{i\in\ZZ})_j=a_{j+1}$ for morphisms.

We embed $\calE$ into $\Gr^\mathrm{b}(\calE)$ by sending objects
(and morphisms) to tuples concentrated in degree $0$. If
$(g_1,g_2,\mu)$ is a determinant functor in
$\det\big((\Gr^\mathrm{b}(\calE),T),\calP\big)$ then the restriction
of $(g_1,g_2)$ to $\calE$ is a determinant functor in
$\det(\calE,\calP)$. Clearly we can extend this to a functor $W:
\det\big((\Gr^\mathrm{b}(\calE),T),\calP\big)\to\det(\calE,\calP)$.

\begin{proposition}
\label{proposition_graded} The functor $W:
\det\big((\Gr^\mathrm{b}(\calE),T),\calP\big)\to\det(\calE,\calP)$
is an equivalence of categories.
\end{proposition}

\begin{proof}
If $A=(A_i)$ is an object in $\Gr^\mathrm{b}(\calE)$ then we write
$A^\mathrm{ev}$ and $A^\mathrm{od}$ for the objects
$A^\mathrm{ev}=\bigoplus_{i\,\mathrm{even}}A_i$ and
$A^\mathrm{od}=\bigoplus_{i\,\mathrm{odd}}A_i$ in $\calE$. For a
morphism $a=(a_i): A\to B$ in $\Gr^\mathrm{b}(\calE)$ we obtain
induced morphisms $a^\mathrm{ev}: A^\mathrm{ev}\to B^\mathrm{ev}$
and $a^\mathrm{od}: A^\mathrm{od}\to B^\mathrm{od}$. This defines
two exact functors $\Gr^\mathrm{b}(\calE)\to\calE$.

For every object $X$ in $\calP$ we fix a right inverse $X^*$, i.e.\
an object $X^*$ together with a unit structure on $X\otimes X^*$. If
$\alpha: X\to Y$ is a morphism in $\calP$ then we let $\alpha^*:
X^*\to Y^*$ be the unique morphism for which $\alpha\otimes\alpha^*:
X\otimes X^*\to Y\otimes Y^*$ is a morphism of units. Furthermore we
note that there exists a canonical isomorphism $(X\otimes Y)^*\cong
X^*\otimes Y^*$.

If $(f_1,f_2)$ is a determinant functor in $\det(\calE,\calP)$ then
we define a functor $g_1: \Gr^\mathrm{b}(\calE)_\iso\to\calP$ by
$g_1(A)=f_1(A^\mathrm{ev})\otimes f_1(A^\mathrm{od})^*$ for objects
$A$ in $\Gr^\rmb(\calE)$ and similarly for isomorphisms. For every
short exact sequence $\Delta: A\to B\to C$ in $\Gr^\rmb(\calE)$ we
define $g_2(\Delta): g_1(B)\to g_1(A)\otimes g_1(C)$ to be the
composite morphism
\begin{equation*}
\begin{split}
g_1(B) & = f_1(B^\mathrm{ev})\otimes f_1(B^\mathrm{od})^* \\
&  \xrightarrow{f_2(\Delta^\mathrm{ev})\otimes
f_2(\Delta^\mathrm{od})^*} \big(f_1(A^\mathrm{ev})\otimes
f_1(C^\mathrm{ev})\big)\otimes\big(f_1(A^\mathrm{od})\otimes
f_1(C^\mathrm{od})\big)^* \\
& \xrightarrow{\,\cong\,} \big(f_1(A^\mathrm{ev})\otimes
f_1(C^\mathrm{ev})\big)\otimes\big(f_1(A^\mathrm{od})^*\otimes
f_1(C^\mathrm{od})^*\big) \\
& \xrightarrow{\,\cong\,} \big(f_1(A^\mathrm{ev})\otimes
f_1(A^\mathrm{od})^*\big)\otimes\big(f_1(C^\mathrm{ev})\otimes
f_1(C^\mathrm{od})^*\big) = g_1(A)\otimes g_1(C),
\end{split}
\end{equation*}
where $\Delta^\mathrm{ev}: A^\mathrm{ev}\to B^\mathrm{ev}\to
C^\mathrm{ev}$ and $\Delta^\mathrm{od}: A^\mathrm{od}\to
B^\mathrm{od}\to C^\mathrm{od}$. Finally we define $\mu_A$ to be the
unit structure on $g_1(A)\otimes g_1(TA)$ which is induced from the
product unit structure on $\big(f_1(A^\mathrm{ev})\otimes
f_1(A^\mathrm{ev})^*\big)\otimes \big(f_1(A^\mathrm{od})\otimes
f_1(A^\mathrm{od})^*\big)$ via the isomorphism
\begin{equation*}
\begin{split}
g_1(A)\otimes g_1(TA) & = \big(f_1(A^\mathrm{ev})\otimes
f_1(A^\mathrm{od})^*\big)\otimes \big(f_1(A^\mathrm{od})\otimes
f_1(A^\mathrm{ev})^*\big) \\
& \xrightarrow{\,\cong\,} \big(f_1(A^\mathrm{ev})\otimes
f_1(A^\mathrm{ev})^*\big)\otimes \big(f_1(A^\mathrm{od})\otimes
f_1(A^\mathrm{od})^*\big) \\
& \xrightarrow{\varepsilon(f_1(A^\mathrm{od}))}
\big(f_1(A^\mathrm{ev})\otimes f_1(A^\mathrm{ev})^*\big)\otimes
\big(f_1(A^\mathrm{od})\otimes f_1(A^\mathrm{od})^*\big).
\end{split}
\end{equation*}
One can check that $(g_1,g_2,\mu)$ is a determinant functor in
$\det\big((\Gr^\mathrm{b}(\calE),T),\calP\big)$. Furthermore
$(f_1,f_2)\mapsto (g_1,g_2,\mu)$ naturally extends to a functor $R:
\det(\calE,\calP)\to\det\big((\Gr^\mathrm{b}(\calE),T),\calP\big)$.

It is straightforward to verify that the composite functor
\begin{equation*}
\det(\calE,\calP) \xrightarrow{R}
\det\big((\Gr^\mathrm{b}(\calE),T),\calP\big) \xrightarrow{W}
\det(\calE,\calP)
\end{equation*}
is isomorphic to the identity functor on $\det(\calE,\calP)$.

Next we want to show that the composite of the functors
\begin{equation*}
\det\big((\Gr^\mathrm{b}(\calE),T),\calP\big) \xrightarrow{W}
\det(\calE,\calP) \xrightarrow{R}
\det\big((\Gr^\mathrm{b}(\calE),T),\calP\big)
\end{equation*}
is isomorphic to the identity functor on
$\det\big((\Gr^\mathrm{b}(\calE),T),\calP\big)$. For this we use the
results from \cite[\S 4]{BreuningBurns05}. Note that those results
are formulated for determinant functors on certain categories of
cochain complexes, but that they also apply to determinant functors
on $(\Gr^\mathrm{b}(\calE),T)$. If $(g_1,g_2,\mu)$ is in
$\det\big((\Gr^\mathrm{b}(\calE),T),\calP\big)$ and
$R(W(g_1,g_2,\mu))=(g'_1,g'_2,\mu')$, then from \cite[Prop.\
4.4]{BreuningBurns05} we obtain an isomorphism of determinant
functors $\pi: (g_1,g_2)\to (g'_1,g'_2)$ in
$\det(\Gr^\rmb(\calE),\calP)$. Furthermore \cite[Lemma
4.5]{BreuningBurns05} shows that $\pi$ is compatible with $\mu$ and
$\mu'$ and is therefore a morphism of determinant functors in
$\det\big((\Gr^\rmb(\calE),T),\calP\big)$. Finally it is not
difficult to verify that $\pi$ is natural in $g$ and hence gives a
natural isomorphism from the identity functor on
$\det\big((\Gr^\rmb(\calE),T),\calP\big)$ to $R\circ W$
\end{proof}

\subsection{Restriction functor to the graded category}
\label{subsection_graded_restriction}

Let $\calT$ be a triangulated category with a t-structure and let
$\calE$ denote its heart. Let $f=(f_1,f_2): \calT\to\calP$ be a
determinant functor on $\calT$. We want to show that $f$ induces a
determinant functor $V(f)$ on the exact category with translation
$(\Gr^\mathrm{b}(\calE),T)$.

First we define a functor $J:\Gr^\mathrm{b}(\calE)\to\calT$ as
follows. If $A=(A_i)_{i\in\ZZ}$ is an object in
$\Gr^\mathrm{b}(\calE)$ then $J(A)=\bigoplus_{i\in\ZZ}T^{-i}A_i$ and
if $a=(a_i)_{i\in\ZZ}: A\to B$ is a morphism in
$\Gr^\mathrm{b}(\calE)$ then $J(a)=\bigoplus_{i\in\ZZ}T^{-i}a_i:
J(A)\to J(B)$. Note that there exists a canonical natural
isomorphism $TJ\cong JT$. Let $\Delta: 0\to A\xrightarrow{a}
B\xrightarrow{b} C\to 0$ be a short exact sequence in
$\Gr^\mathrm{b}(\calE)$. Then for every $i\in\ZZ$ we have a short
exact sequence $0\to A_i\xrightarrow{a_i} B_i\xrightarrow{b_i}
C_i\to 0$ in $\calE$ and hence there is a unique morphism $c_i:
C_i\to TA_i$ in $\calT$ such that $A_i\xrightarrow{a_i}
B_i\xrightarrow{b_i} C_i\xrightarrow{c_i} TA_i$ is a distinguished
triangle. We define $J(\Delta)$ to be the distinguished triangle
$J(A)\xrightarrow{J(a)} J(B)\xrightarrow{J(b)} J(C)\xrightarrow{c'}
TJ(A)$ where $c'=\bigoplus_{i\in\ZZ} (-1)^i T^{-i}c_i$.

We define a functor $g_1: \Gr^\mathrm{b}(\calE)_\iso\to\calP$ by
$g_1=f_1\circ J$, and for every short exact sequence $\Delta: 0\to
A\to B\to C\to 0$ in $\Gr^\mathrm{b}(\calE)$ we define an
isomorphism $g_2(\Delta): g_1(B)\to g_1(A)\otimes g_1(C)$ by
$g_2(\Delta)=f_2(J(\Delta))$. For every object $A$ in
$\Gr^\mathrm{b}(\calE)$ we define a unit structure $\mu'_A$ on
$g_1(A)\otimes g_1(TA)$ by $\mu'_A=\mu_{J(A)}$ where $\mu_{J(A)}$ is
the unit structure on $f_1(J(A))\otimes f_1(TJ(A))$ defined in Lemma
\ref{lemma_mu}. Let $V(f)$ denote the triple $(g_1,g_2,\mu')$.

\begin{lemma}
If $f$ is a determinant functor on $\calT$ with values in $\calP$,
then $V(f)$ is a determinant functor on $(\Gr^\mathrm{b}(\calE),T)$
with values in $\calP$. Moreover $f\mapsto V(f)$ naturally extends
to a functor $V: \det(\calT,\calP)\to
\det\big((\Gr^\mathrm{b}(\calE),T),\calP\big)$.
\end{lemma}

\begin{proof}
To see that $(g_1,g_2)$ is a determinant functor on the exact
category $\Gr^\mathrm{b}(\calE)$ one proceeds as in the proof of
Lemma \ref{lemma_restriction}. The necessary properties of the unit
structure $\mu'$ follow from Lemma \ref{lemma_mu}. The extension of
$f\mapsto V(f)$ to a functor $V: \det(\calT,\calP)\to
\det\big((\Gr^\mathrm{b}(\calE),T),\calP\big)$ is again similar to
the proof of Lemma \ref{lemma_restriction}.
\end{proof}

\begin{lemma}
The composite functor $W\circ V:
\det(\calT,\calP)\to\det(\calE,\calP)$ is isomorphic to the
restriction functor $f\mapsto f|_{\calE}$.
\end{lemma}

\begin{proof}
This is immediate from the definitions.
\end{proof}

\subsection{Cohomological functors}
\label{subsection_cohomological_functor}

Let $\calT$ be a triangulated category, $\calE$ an abelian category
and $H: \calT\to\calE$ a cohomological functor, i.e.\ an additive
functor such that for every distinguished triangle $\Delta:
A\xrightarrow{a} B\xrightarrow{b} C\xrightarrow{c} TA$ in $\calT$
the long sequence
\begin{equation}
\label{eqn_long_cohom_sequence}
\begin{split}
\dots \xrightarrow{H(T^{i-1}c)} & H(T^iA) \xrightarrow{H(T^ia)}
H(T^iB) \xrightarrow{H(T^ib)} H(T^iC) \\
\xrightarrow{H(T^ic)} & H(T^{i+1}A) \xrightarrow{H(T^{i+1}a)} \dots
\end{split}
\end{equation}
is exact in $\calE$. We assume that $H$ is bounded, i.e.\ that for
every object $A$ of $\calT$ we have $H(T^iA)=0$ for all but finitely
many $i\in\ZZ$. Then we can define a functor $\tilde{H}:
\calT\to\Gr^\mathrm{b}(\calE)$ by setting
$\tilde{H}(A)=(H(T^iA))_{i\in\ZZ}$ (resp.\
$\tilde{H}(a)=(H(T^ia))_{i\in\ZZ}$) for objects $A$ (resp.\
morphisms $a$) in $\calT$. The long exact sequence
(\ref{eqn_long_cohom_sequence}) can be written as the exact sequence
\begin{equation}
\label{eqn_short_cohom_sequence} 0\to \ker\tilde{H}a \to \tilde{H}A
\xrightarrow{\tilde{H}a} \tilde{H}B \xrightarrow{\tilde{H}b}
\tilde{H}C \xrightarrow{\tilde{H}c} T\ker\tilde{H}a \to 0
\end{equation}
in $\Gr^\mathrm{b}(\calE)$.

Now consider a determinant functor $g=(g_1,g_2,\mu):
(\Gr^\mathrm{b}(\calE),T)\to\calP$ on the exact category with
translation $(\Gr^\mathrm{b}(\calE),T)$. We define a functor $f_1:
\calT_\iso\to\calP$ as the composite $f_1=g_1\circ \tilde{H}$. Let
$\Delta: A\xrightarrow{a} B\xrightarrow{b} C\xrightarrow{c} TA$ be a
distinguished triangle in $\calT$ and write the long exact
cohomology sequence of $\Delta$ as the exact sequence
(\ref{eqn_short_cohom_sequence}) in $\Gr^\mathrm{b}(\calE)$. We
define $f_2(\Delta): f_1(B)\to f_1(A)\otimes f_1(C)$ to be the
composite
\begin{equation*}
\begin{split}
f_1(B) = g_1(\tilde{H}B) & \xrightarrow{\mu_{\ker\tilde{H}a}}
g_1(\tilde{H}B)\otimes g_1(\ker\tilde{H}a)\otimes g_1(T\ker\tilde{H}a) \\
& \xrightarrow{g_2(\Delta_1)\otimes\id\otimes\id}
g_1(\ker\tilde{H}b)\otimes g_1(\ker\tilde{H}c)\otimes
g_1(\ker\tilde{H}a)\otimes g_1(T\ker\tilde{H}a) \\
& \xrightarrow{\,\cong\,} g_1(\ker\tilde{H}a)\otimes
g_1(\ker\tilde{H}b)\otimes g_1(\ker\tilde{H}c)\otimes
g_1(T\ker\tilde{H}a) \\
& \xrightarrow{g_2(\Delta_2)^{-1}\otimes g_2(\Delta_3)^{-1}}
g_1(\tilde{H}A)\otimes g_1(\tilde{H}C) = f_1(A)\otimes f_1(C),
\end{split}
\end{equation*}
where the short exact sequence $\Delta_1: 0\to\ker\tilde{H}b\to
\tilde{H}B\to \ker\tilde{H}c\to 0$, $\Delta_2: 0\to
\ker\tilde{H}a\to \tilde{H}A\to \ker\tilde{H}b\to 0$ and $\Delta_3:
0\to \ker\tilde{H}c\to \tilde{H}C\to T\ker\tilde{H}a\to 0$ are
obtained by splitting (\ref{eqn_short_cohom_sequence}). Let $H^*(g)$
denote the pair $(f_1,f_2)$.

\begin{proposition}
\label{proposition_cohom_functor} If $g$ is a determinant functor in
$\det\big((\Gr^\mathrm{b}(\calE),T),\calP\big)$ then $H^*(g)$ is a
determinant functor in $\det(\calT,\calP)$. Moreover $g\mapsto
H^*(g)$ naturally extends to a functor $H^*:
\det\big((\Gr^\mathrm{b}(\calE),T),\calP\big)\to\det(\calT,\calP)$.
\end{proposition}

\begin{remark}
The definition of $H^*(g)$ and Proposition
\ref{proposition_cohom_functor} are also valid for an unbounded
cohomological functor $H$ if we use the category $\Gr(\calE)$ of all
$\ZZ$-graded objects in $\calE$ instead of $\Gr^\mathrm{b}(\calE)$.
\end{remark}

\begin{proof}[Proof of Proposition \ref{proposition_cohom_functor}]
The second statement of the proposition follows easily once we have
proved the first statement. Indeed, let $g=(g_1,g_2,\mu)$ and
$g'=(g'_1,g'_2,\mu')$ be in
$\det\big((\Gr^\mathrm{b}(\calE),T),\calP\big)$ and define
$H^*(g)=(f_1,f_2)$ and $H^*(g')=(f'_1,f'_2)$ as above. If $\lambda:
g\to g'$ is a morphism in
$\det\big((\Gr^\mathrm{b}(\calE),T),\calP\big)$ then for every
object $A$ in $\calT$ we can define
$H^*(\lambda)_A=\lambda_{\tilde{H}A}: f_1(A)=g_1(\tilde{H}A)\to
g'_1(\tilde{H}A)=f'_1(A)$. It is straightforward to check that
$H^*(\lambda)$ is a natural isomorphism $f_1\to f'_1$ and that
$(H^*(\lambda)_A\otimes H^*(\lambda)_C)\circ f_2(\Delta) =
f'_2(\Delta)\circ H^*(\lambda)_B$ for every distinguished triangle
$\Delta: A\to B\to C\to TA$ in $\calT$. Hence if we know that
$H^*(g)$ and $H^*(g')$ are determinant functors, then $H^*(\lambda):
H^*(g)\to H^*(g')$ is a morphism of determinant functors, and we
therefore obtain a functor $H^*:
\det\big((\Gr^\mathrm{b}(\calE),T),\calP\big)\to\det(\calT,\calP)$.

The proof of the first statement of Proposition
\ref{proposition_cohom_functor} will occupy the remainder of \S
\ref{subsection_cohomological_functor}. We let $g=(g_1,g_2,\mu)$ be
a determinant functor in
$\det\big((\Gr^\mathrm{b}(\calE),T),\calP\big)$ and define
$H^*(g)=(f_1,f_2)$ as above. To simplify the notation we will write
$H$ instead of $\tilde{H}$ for the rest of this proof.

An isomorphism of distinguished triangles induces an isomorphism of
long exact cohomology sequences. From this it follows that $f_2$ is
natural in isomorphisms of distinguished triangles.

The cohomology sequence (\ref{eqn_short_cohom_sequence}) associated
to a distinguished triangle $\Delta$ of the form $A\to A\oplus B\to
B\xrightarrow{0} TA$ is of the form $0\to 0\to HA\to H(A\oplus B)\to
HB\to T0\to 0$. It follows easily that $f_2(\Delta)=g_2(\Delta')$
where $\Delta'$ is the short exact sequence $HA\to H(A\oplus B)\to
HB$ in $\Gr^\mathrm{b}(\calE)$. Since $H(A\oplus B)$ is a direct sum
of $HA$ and $HB$ and the maps $HA\to H(A\oplus B)$ and $H(A\oplus
B)\to HB$ are the canonical inclusion and projection, we can deduce
the commutativity axiom for $H^*(g)$ from the commutativity axiom
for $g$.

It remains to show the associativity axiom. So assume that we are
given the following octahedral diagram in $\calT$.
\begin{equation}
\label{diagram_octahedral_cohom}
\begin{split}
\xymatrix{ A \ar[r]^a \ar@{=}[d] & B \ar[r]^{b} \ar[d]^{c} & C'
\ar[r]^{fd} \ar[d]^{d} & TA \ar@{=}[d] \\
A \ar[r]^{ca} & C \ar[r]^{e} \ar[d]^{ge} & B' \ar[r]^{f} \ar[d]^{g}
& TA \\
& A' \ar@{=}[r] \ar[d]^{h} & A' \ar[d]^{Tb\circ h} & \\
& TB \ar[r]^{Tb} & TC' & }
\end{split}
\end{equation}
We must show that
\begin{equation*}
\xymatrix@C+1.5cm{ f_1(C) \ar[r]^-{f_2(\Delta_{\mathrm{h2}})}
\ar[d]^{f_2(\Delta_{\mathrm{v1}})} & f_1(A)\otimes f_1(B')
\ar[d]^{\id\otimes f_2(\Delta_{\mathrm{v2}})} \\
f_1(B)\otimes f_1(A') \ar[r]^-{f_2(\Delta_{\mathrm{h1}})\otimes\id}
& f_1(A)\otimes f_1(C')\otimes f_1(A') }
\end{equation*}
is commutative in $\calP$, where $\Delta_{\mathrm{h1}}$ and
$\Delta_{\mathrm{h2}}$ (resp.\ $\Delta_{\mathrm{v1}}$ and
$\Delta_{\mathrm{v2}}$) are the first and second horizontal (resp.\
vertical) distinguished triangles in
(\ref{diagram_octahedral_cohom}). By definition of $f_1$ and $f_2$
this is equivalent to showing that
\begin{equation}
\label{diagram_cohom_assoc}
\begin{split}
\xymatrix@C+2.2cm@R+1.2cm{ g_1(HC) \ar[r]^-{\begin{array}{c}
\scriptstyle \mu_{\ker Hca}, g_2(\Delta_1),\\ \scriptstyle
g_2(\Delta_2)^{-1}, g_2(\Delta_3)^{-1} \end{array}}
\ar[d]^{\begin{array}{c} \scriptstyle \mu_{\ker Hc},
g_2(\Delta_4),\\ \scriptstyle g_2(\Delta_5)^{-1}, g_2(\Delta_6)^{-1}
\end{array}} & g_1(HA)\otimes g_1(HB')
\ar[d]^{\begin{array}{c} \scriptstyle \mu_{\ker Hd}, g_2(\Delta_7),\\
\scriptstyle g_2(\Delta_8)^{-1},
g_2(\Delta_9)^{-1} \end{array}} \\
g_1(HB)\otimes g_1(HA') \ar[r]^-{\begin{array}{c} \scriptstyle
\mu_{\ker Ha}, g_2(\Delta_{10}),\\ \scriptstyle
g_2(\Delta_{11})^{-1}, g_2(\Delta_{12})^{-1}
\end{array}} & g_1(HA)\otimes g_1(HC')\otimes g_1(HA') }
\end{split}
\end{equation}
commutes in $\calP$, where the short exact sequences $\Delta_1,
\dots, \Delta_{12}$ are obtained by splitting the exact cohomology
sequences associated to the distinguished triangles in
(\ref{diagram_octahedral_cohom}). To show the commutativity of
(\ref{diagram_cohom_assoc}) we must construct a huge commutative
diagram in $\calP$, which we do in several pieces.

First we observe that applying $H$ to diagram
(\ref{diagram_octahedral_cohom}) gives the following commutative
diagram with exact rows and columns in $\Gr^\rmb(\calE)$
\begin{equation*}
\xymatrix{
& & & 0 \ar[d] & 0 \ar[d] & & \\
& & & \ker Hc \ar[r] \ar[d] & \ker Hd \ar[d] & & \\
0 \ar[r] & \ker Ha \ar[r] \ar[d] & HA \ar[r]^{Ha} \ar@{=}[d] & HB
\ar[r]^{Hb} \ar[d]^{Hc} & HC' \ar[r]^{Hfd} \ar[d]^{Hd} & T\ker Ha
\ar[r] \ar[d] & 0 \\
0 \ar[r] & \ker Hca \ar[r] & HA \ar[r]^{Hca} & HC \ar[r]^{He}
\ar[d]^{Hge} & HB' \ar[r]^{Hf} \ar[d]^{Hg} & T\ker Hca \ar[r] & 0 \\
& & & HA' \ar@{=}[r] \ar[d]^{Hh} & HA' \ar[d]^{H(Tb\circ h)} & & \\
& & & T\ker Hc \ar[r]^{HTb} \ar[d] & T\ker Hd \ar[d] & & \\
& & & 0 & 0 & & }
\end{equation*}
and that applying $H$ to the equality $Ta\circ f=h\circ g$ gives
\begin{equation}
\label{eqn_special_square}
THa\circ Hf=Hh\circ Hg.
\end{equation}
In $\Gr^\rmb(\calE)$ we have the following short exact sequences.
\begin{gather*}
\Delta_{1}: \ker He\to HC\xrightarrow{He} \ker Hf \\
\Delta_{2}: \ker Hca\to HA\xrightarrow{Hca} \ker He \\
\Delta_{3}: \ker Hf\to HB'\xrightarrow{Hf} T\ker Hca \\
\Delta_{4}: \ker Hge\to HC\xrightarrow{Hge} \ker Hh \\
\Delta_{5}: \ker Hc\to HB\xrightarrow{Hc} \ker Hge \\
\Delta_{6}: \ker Hh\to HA'\xrightarrow{Hh} T\ker Hc \\
\Delta_{7}: \ker Hg\to HB'\xrightarrow{Hg} \ker H(Tb\circ h) \\
\Delta_{8}: \ker Hd\to HC'\xrightarrow{Hd} \ker Hg \\
\Delta_{9}: \ker H(Tb\circ h)\to HA'\xrightarrow{H(Tb\circ h)} T\ker
Hd
\end{gather*}
\begin{gather*}
\Delta_{10}: \ker Hb\to HB\xrightarrow{Hb} \ker Hfd \\
\Delta_{11}: \ker Ha\to HA\xrightarrow{Ha} \ker Hb \\
\Delta_{12}: \ker Hfd\to HC'\xrightarrow{Hfd} T\ker Ha \\
\Delta_{13}: \ker Ha\to \ker Hca\xrightarrow{Ha} \ker Hb\cap\ker Hc \\
\Delta_{14}: \ker Hb\cap\ker Hc\to \ker Hc\xrightarrow{Hb} \ker Hd
\end{gather*}
Furthermore we will need the translations of $\Delta_{13}$ and
$\Delta_{14}$, that is the following two short exact sequences.
\begin{gather*}
T\Delta_{13}: T\ker Ha\to T\ker Hca\xrightarrow{THa} T\ker Hb\cap T\ker Hc \\
T\Delta_{14}: T\ker Hb\cap T\ker Hc\to T\ker Hc\xrightarrow{THb}
T\ker Hd
\end{gather*}

We will show the commutativity of the following diagram. Since its
outer square is (\ref{diagram_cohom_assoc}) this will imply the
associativity axiom for $H^*(g)$. For typographic reasons we will
write $[\cdot]$ for $g_1(\cdot)$ and omit the word $\ker$ in all
diagrams below, so for example we have $[HA]=g_1(HA)$ and
$[Ha]=g_1(\ker Ha)$. Furthermore we will write $\Delta_i$ for both
$g_2(\Delta_i)$ and $g_2(\Delta_i)^{-1}$.
\begin{equation*} \xymatrix@C+1.1cm{ \scriptstyle [HC]
\ar[r]^-{\mu_{Hca}} \ar[d]^{\mu_{Hc}} \ar[dr]^{\mu_{Hd}, \mu_{Ha},
\mu_{Hb\cap Hc}} & {\begin{array}{c} \scriptstyle [HC] \\
\scriptstyle \otimes [Hca]\otimes [THca] \end{array}}
\ar[r]^-{\Delta_1, \Delta_2, \Delta_3} \ar[d]^{\mu_{Hd},
\Delta_{13}, T\Delta_{13}} &
\scriptstyle [HA]\otimes [HB'] \ar[d]^{\mu_{Hd}} \\
{\begin{array}{c} \scriptstyle [HC] \\ \scriptstyle \otimes
[Hc]\otimes [THc] \end{array}} \ar[r]^-{\mu_{Ha}, \Delta_{14},
T\Delta_{14}} \ar[d]^{\Delta_4, \Delta_5, \Delta_6} &
{\begin{array}{c} \scriptstyle [HC] \\
\scriptstyle \otimes [Ha]\otimes [THa] \\
\scriptstyle \otimes [Hb\cap Hc]\otimes [THb\cap THc]
\\ \scriptstyle \otimes [Hd]\otimes [THd]
\end{array}} \ar[r]^-{\begin{array}{c} \scriptstyle \Delta_{13}, T\Delta_{13},\\
\scriptstyle \Delta_1, \Delta_2, \Delta_3 \end{array}}
\ar[d]^{\Delta_{14}, T\Delta_{14}, \Delta_4, \Delta_5, \Delta_6} &
{\begin{array}{c} \scriptstyle [HA]\otimes [HB'] \\
\scriptstyle \otimes [Hd]\otimes [THd] \end{array}}
\ar[d]^{\Delta_7, \Delta_8, \Delta_9} \\
\scriptstyle [HB]\otimes [HA'] \ar[r]^-{\mu_{Ha}} &
{\begin{array}{c} \scriptstyle
[HB]\otimes [HA'] \\
\scriptstyle \otimes [Ha]\otimes [THa] \end{array}}
\ar[r]^-{\Delta_{10}, \Delta_{11}, \Delta_{12}} & \scriptstyle
[HA]\otimes [HC']\otimes [HA'] }
\end{equation*}
From the properties of $\mu$ it immediately follows that the two
triangles are commutative. Also the commutativity of the top right
and bottom left hand squares is obvious.

It remains to show that the bottom right hand square is commutative.
This square can be subdivided as follows.
\begin{equation*}
\xymatrix@C+0.7cm@R+0.2cm{ {\begin{array}{c}
\scriptstyle [HC] \\
\scriptstyle \otimes [Ha]\otimes [THa] \\
\scriptstyle \otimes [Hb\cap Hc]\otimes [THb\cap THc]
\\ \scriptstyle \otimes [Hd]\otimes [THd]
\end{array}} \ar[r]^-{\Delta_1} \ar[d]^{\Delta_4}
\ar@{}[dr]|{\textstyle S_1} &
{\begin{array}{c} \scriptstyle [He]\otimes [Hf] \\
\scriptstyle \otimes [Ha]\otimes [THa] \\
\scriptstyle \otimes [Hb\cap Hc]\otimes [THb\cap THc]
\\ \scriptstyle \otimes [Hd]\otimes [THd]
\end{array}} \ar[r]^-{\begin{array}{c} \scriptstyle \Delta_{13}, T\Delta_{13},\\
\scriptstyle \Delta_2, \Delta_3 \end{array}} \ar[d]^{\Delta_{15}}
\ar@{}[dr]|{\textstyle
S_2} & {\begin{array}{c} \scriptstyle [HA]\otimes [HB'] \\
\scriptstyle \otimes [Hd]\otimes
[THd] \end{array}} \ar[d]^{\Delta_7} \\
{\begin{array}{c} \scriptstyle [Hge]\otimes [Hh] \\
\scriptstyle \otimes [Ha]\otimes [THa] \\
\scriptstyle \otimes [Hb\cap Hc]\otimes [THb\cap THc]
\\ \scriptstyle \otimes [Hd]\otimes [THd]
\end{array}}
\ar[r]^-{\Delta_{16}} \ar[d]^-{\begin{array}{c} \scriptstyle
\Delta_{14}, T\Delta_{14},\\ \scriptstyle \Delta_5, \Delta_6
\end{array}} \ar@{}[dr]|{\textstyle S_3} & {\begin{array}{c}
\scriptstyle [He]\otimes [Hf\cap Hg]\otimes [Hh] \\
\scriptstyle \otimes [Ha]\otimes [THa] \\
\scriptstyle \otimes [Hb\cap Hc]\otimes [THb\cap THc]
\\ \scriptstyle \otimes [Hd]\otimes [THd]
\end{array}} \ar[r]^-{\begin{array}{c} \scriptstyle \Delta_2, \Delta_{13},\\
\scriptstyle \Delta_{17}, \Delta_{18} \end{array}}
\ar[d]^-{\begin{array}{c} \scriptstyle \Delta_6, T\Delta_{14},\\
\scriptstyle \Delta_{19}, \Delta_{20} \end{array}}
\ar@{}[dr]|{\textstyle S_4}
& {\begin{array}{c} \scriptstyle [HA] \\
\scriptstyle \otimes [Hg]\otimes [H(Tb\circ h)] \\
\scriptstyle \otimes [Hd]\otimes [THd] \end{array}}
\ar[d]^{\Delta_8, \Delta_9} \\
{\begin{array}{c} \scriptstyle [HB]\otimes [HA'] \\
\scriptstyle \otimes [Ha]\otimes [THa] \end{array}}
\ar[r]^-{\Delta_{10}} & {\begin{array}{c} \scriptstyle [Hb]\otimes
[Hfd]
\\ \scriptstyle \otimes [Ha]\otimes [THa] \\ \scriptstyle \otimes [HA']
\end{array}} \ar[r]^-{\Delta_{11}, \Delta_{12}} &
\scriptstyle [HA]\otimes [HC']\otimes [HA'] }
\end{equation*}
Here $\Delta_{15}, \dots, \Delta_{20}$ are the following short exact
sequences in $\Gr^\rmb(\calE)$.
\begin{gather*}
\Delta_{15}: \ker Hf\cap\ker Hg\to \ker Hf\xrightarrow{Hg} \ker Hh
\\
\Delta_{16}: \ker He\to \ker Hge\xrightarrow{He} \ker Hf\cap\ker Hg
\\
\Delta_{17}: \ker Hf\cap\ker Hg\to \ker Hg\xrightarrow{Hf} T\ker Ha \\
\Delta_{18}: \ker Hh\to \ker H(Tb\circ h)\xrightarrow{Hh} T\ker
Hb\cap T\ker Hc \\
\Delta_{19}: \ker Hb\cap\ker Hc\to \ker Hb\xrightarrow{Hc} \ker He \\
\Delta_{20}: \ker Hd\to \ker Hfd\xrightarrow{Hd} \ker Hf\cap\ker Hg
\end{gather*}

Applying the associativity axiom to the commutative diagram of short
exact sequences
\begin{equation*}
\xymatrix{ \ker He \ar[r] \ar@{=}[d] & \ker Hge \ar[r]^-{He} \ar[d]
& \ker Hf\cap\ker Hg \ar[d] \\
\ker He \ar[r] & HC \ar[r]^-{He} \ar[d]^{Hge} & \ker Hf \ar[d]^{Hg}
\\
& \ker Hh \ar@{=}[r] & \ker Hh }
\end{equation*}
gives the commutativity of the square $S_1$.

Next we will show the commutativity of the square $S_2$. For this we
subdivide $S_2$ as in the following diagram
\begin{equation*}
\xymatrix@C+0.8cm{ {\begin{array}{c} \scriptstyle [He]\otimes [Hf] \\
\scriptstyle \otimes [Ha]\otimes [THa] \\
\scriptstyle \otimes [Hb\cap Hc]\otimes [THb\cap THc]
\\ \scriptstyle \otimes [Hd]\otimes [THd]
\end{array}} \ar[rr]^-{\Delta_{13}, \Delta_2, T\Delta_{13}, \Delta_3}
\ar[dr]^{\Delta_{13}, \Delta_2, \Delta_{21}} \ar[dd]^{\Delta_{15}} &
\ar@{}[d]|-{\textstyle S_{2,1}} &
{\begin{array}{c} \scriptstyle [HA]\otimes [HB'] \\
\scriptstyle \otimes [Hd]\otimes [THd]
\end{array}} \ar[dd]^{\Delta_7} \ar[dl]^{\Delta_{22}} \\
\ar@{}[r]_-{\textstyle S_{2,3}} & {\begin{array}{c} \scriptstyle [HA]\otimes [H(Ta\circ f)] \\
\scriptstyle \otimes [THb\cap THc] \\
\scriptstyle \otimes [Hd] \otimes [THd] \end{array}} \ar[d]^{\Delta_{23}}
\ar@{}[r]_-{\textstyle S_{2,2}} & \\
{\begin{array}{c} \scriptstyle [He] \otimes [Hf\cap Hg]\otimes
[Hh] \\ \scriptstyle \otimes [Ha]\otimes [THa] \\
\scriptstyle \otimes
[Hb\cap Hc]\otimes [THb\cap THc] \\
\scriptstyle \otimes [Hd]\otimes [THd] \end{array}}
\ar[r]^-{\Delta_{13}, \Delta_2, \Delta_{17}} &
{\begin{array}{c} \scriptstyle [HA]\otimes [Hg]\otimes [Hh] \\
\scriptstyle \otimes [THb\cap THc] \\
\scriptstyle \otimes [Hd]\otimes [THd] \end{array}}
\ar[r]^-{\Delta_{18}} & {\begin{array}{c} \scriptstyle [HA]\otimes
[Hg]\otimes [H(Tb\circ h)] \\ \scriptstyle \otimes [Hd]\otimes [THd]
\end{array}} }
\end{equation*}
where $\Delta_{21}, \dots, \Delta_{23}$ are the short exact
sequences
\begin{gather*}
\Delta_{21}: \ker Hf\to \ker H(Ta\circ f)\xrightarrow{Hf} T\ker Ha
\\
\Delta_{22}: \ker H(Ta\circ f)\to HB'\xrightarrow{H(Ta\circ f)}
T\ker Hb\cap T\ker Hc \\
\Delta_{23}: \ker Hg\to \ker H(Ta\circ f)\xrightarrow{Hg} \ker Hh.
\end{gather*}
The commutativity of $S_{2,1}$ and $S_{2,2}$ comes from the
associativity axiom applied to the following two commutative
diagrams of short exact sequences
\begin{equation*}
\xymatrix{ \ker Hf \ar[r] \ar@{=}[d] & \ker H(Ta\circ f) \ar[r]^{Hf}
\ar[d] &
T\ker Ha \ar[d] \\
\ker Hf \ar[r] & HB' \ar[r]^{Hf} \ar[d]^{H(Ta\circ f)} & T\ker Hca
\ar[d]^{THa} \\
& T\ker Hb\cap T\ker Hc \ar@{=}[r] & T\ker Hb\cap T\ker Hc, }
\end{equation*}
\begin{equation*}
\xymatrix{ \ker Hg \ar[r] \ar@{=}[d] & \ker H(Ta\circ f) \ar[r]^{Hg}
\ar[d] &
\ker Hh \ar[d] \\
\ker Hg \ar[r] & HB' \ar[r]^{Hg} \ar[d]^{H(Ta\circ f)} & \ker
H(Tb\circ h) \ar[d]^{Hh} \\
& T\ker Hb\cap T\ker Hc \ar@{=}[r] & T\ker Hb\cap T\ker Hc, }
\end{equation*}
(where the last square in the second diagram is commutative by
(\ref{eqn_special_square})), and the commutativity of $S_{2,3}$
comes from \cite[Cor.\ 1.10]{Knudsen02} applied to
\begin{equation*}
\xymatrix{ \ker Hf\cap \ker Hg \ar[r] \ar[d] & \ker Hg \ar[r]^-{Hf}
\ar[d] &
T\ker Ha \ar@{=}[d] \\
\ker Hf \ar[r] \ar[d]^{Hg} & \ker H(Ta\circ f) \ar[r]^-{Hf}
\ar[d]^{Hg} & T\ker Ha \\
\ker Hh \ar@{=}[r] & \ker Hh. & }
\end{equation*}

The commutativity of the squares $S_3$ and $S_4$ can be shown
similarly and is left to the reader.
\end{proof}

\subsection{Completion of the proof}
\label{subsection_cohomology}

Let $\calT$ be a triangulated category with a bounded t-structure,
let $\calE$ be its heart and $H:\calT\to\calE$ the canonical
cohomological functor. Recall that in \S
\ref{subsection_graded_restriction} and \S
\ref{subsection_cohomological_functor} we defined functors
\begin{gather*}
V: \det(\calT,\calP) \to \det\big((\Gr^\rmb(\calE),T),\calP\big), \\
H^*: \det\big((\Gr^\rmb(\calE),T),\calP\big) \to \det(\calT,\calP).
\end{gather*}
In this subsection we show that if the t-structure is non-degenerate
and bounded then $V$ is an equivalence of categories with
quasi-inverse $H^*$.

\begin{lemma}
\label{lemma_composite_1} Assume that the t-structure on $\calT$ is
bounded. Then the composite functor
\begin{equation*}
\det\big((\Gr^\rmb(\calE),T),\calP\big) \xrightarrow{H^*}
\det(\calT,\calP) \xrightarrow{V}
\det\big((\Gr^\rmb(\calE),T),\calP\big)
\end{equation*}
is isomorphic to the identity functor on
$\det\big((\Gr^\rmb(\calE),T),\calP\big)$.
\end{lemma}

\begin{proof}
The composite functor $\Gr^\rmb(\calE)\xrightarrow{J}
\calT\xrightarrow{\tilde{H}} \Gr^\rmb(\calE)$ is canonically
isomorphic to the identity functor on $\Gr^\rmb(\calE)$. It follows
that if $g=(g_1,g_2,\mu)\in\det\big((\Gr^\rmb(\calE),T),\calP\big)$
and $V(H^*(g))=(g'_1,g'_2,\mu')$ then there exists a canonical
natural isomorphism $\xi: g_1\to g'_1$. It is easily verified that
$\xi$ is compatible with the isomorphisms $g_2(\Delta)$ and
$g'_2(\Delta)$ for every short exact sequence $\Delta$ in
$\Gr^\rmb(\calE)$ and with the unit structures $\mu_A$ and $\mu'_A$
for every object $A$, and thus is an isomorphism of determinant
functors $\xi: g\to V(H^*(g))$.

It is then straightforward to check that $\xi$ is natural in $g$ and
hence gives a natural isomorphism from the identity functor on
$\det\big((\Gr^\rmb(\calE),T),\calP\big)$ to $V\circ H^*$.
\end{proof}

\begin{proposition}
\label{proposition_cohomology} Assume that the t-structure on
$\calT$ is non-degenerate and bounded. Then the composite functor
\begin{equation*}
\det(\calT,\calP) \xrightarrow{V}
\det\big((\Gr^\rmb(\calE),T),\calP\big) \xrightarrow{H^*}
\det(\calT,\calP)
\end{equation*}
is isomorphic to the identity functor on $\det(\calT,\calP)$.
\end{proposition}

\begin{proof}
We will show that for every determinant functor $f$ in
$\det(\calT,\calP)$ there exists a canonical isomorphism of
determinant functors $\eta: f\to H^*(V(f))$. To complete the proof
of the proposition one must then check that $\eta$ is natural in
$f$. This follows almost immediately from the definition of $\eta$
and is left to the reader.

Write $f=(f_1,f_2)$ and $H^*(V(f))=(f'_1,f'_2)$. We want to
construct a natural isomorphism $\eta: f_1\to f'_1$ that is
compatible with $f_2$ and $f'_2$.

By definition of $H^*$ and $V$ we have $f'_1(A)=f_1(J\tilde{H}A)$
for an object $A$ of $\calT$ and similarly for isomorphisms. Since
the t-structure is non-degenerate and bounded it follows from
\cite[Prop.\ 1.3.7]{BBD} that $A\in\calT^{\geq i}$ for some
$i\in\ZZ$. Then $\tau_{\leq i}A\cong T^{-i}H^iA$ where we write
$H^iA=HT^iA$. Applying $f_2$ to the distinguished triangle
$\tau_{\leq i}A\to A\to\tau_{>i}A\to T\tau_{\leq i}A$ gives an
isomorphism
\begin{equation*}
\begin{split}
f_1(A) & \cong f_1(\tau_{\leq i}A)\otimes f_1(\tau_{>i}A) \\
& \cong f_1(T^{-i}H^iA)\otimes f_1(\tau_{>i}A).
\end{split}
\end{equation*}
Next we apply the same argument to $\tau_{>i}A\in\calT^{\geq i+1}$
and obtain
\begin{equation*}
\begin{split}
f_1(\tau_{>i}A) & \cong
f_1(T^{-(i+1)}H^{i+1}\tau_{>i}A)\otimes f_1(\tau_{>i+1}\tau_{>i}A) \\
& \cong f_1(T^{-(i+1)}H^{i+1}A)\otimes f_1(\tau_{>i+1}A),
\end{split}
\end{equation*}
because $H^{i+1}\tau_{>i}A\cong H^{i+1}A$ and
$\tau_{>i+1}\tau_{>i}A\cong\tau_{>i+1}A$. We then apply the same
argument to $\tau_{>i+1}A\in\calT^{\geq i+2}$, etc. By \cite[Prop.\
1.3.7]{BBD} we have $\tau_{>i+j}A=0$ for sufficiently large $j$, and
therefore we obtain an isomorphism
\begin{equation*}
\begin{split}
\eta_A: f_1(A) & \cong\bigotimes_{i\in\ZZ} f_1(T^{-i}H^iA) \\
& \cong f_1\left(\bigoplus_{i\in\ZZ}T^{-i}H^iA\right) \\
& = f_1(J\tilde{H}A)=f'_1(A).
\end{split}
\end{equation*}
Obviously the isomorphism $\eta_A$ is functorial for isomorphisms
$A\to B$.

To prove that $\eta$ is an isomorphism of determinant functors, we
must show that for every distinguished triangle $\Delta:
A\xrightarrow{a} B\xrightarrow{b} C\xrightarrow{c} TA$ in $\calT$
the diagram
\begin{equation}
\label{diagram_cohomology}
\begin{split}
\xymatrix@C+1cm{ f_1(B) \ar[r]^-{f_2(\Delta)} \ar[d]^{\eta_B} &
f_1(A)\otimes f_1(C) \ar[d]^{\eta_A\otimes\eta_C} \\
f'_1(B) \ar[r]^-{f'_2(\Delta)} & f'_1(A)\otimes f'_1(C) }
\end{split}
\end{equation}
is commutative. Before proving this in full generality we first show
certain special cases and auxiliary results. To simplify the
notation we will write $[A]$ for $f_1(A)$ if $A$ is an object in
$\calT$ and $[A]$ for $f_1(JA)$ if $A$ is an object in
$\Gr^\rmb(\calE)$. In particular we have
$f'_1(A)=f_1(J\tilde{H}A)=[\tilde HA]$.

\begin{lemma}
\label{lemma_cohomology_1} If $A\in\calT^{\geq i}\cap\calT^{\leq i}$
and $B, C\in\calT^{\geq i}$ then diagram (\ref{diagram_cohomology})
is commutative.
\end{lemma}

\begin{proof}
We note that the assumptions imply that the boundary homomorphism in
the long exact cohomology sequence of $\Delta$ (and of all auxiliary
distinguished triangles in this proof) is zero.

Case (i): If $A=\tau_{\leq i}B$, $C=\tau_{>i}B$ then the statement
follows easily from the definition of $\eta_B$.

Case (ii): If $B, C\in\calT^{\leq i}$ then the distinguished
triangle $A\to B\to C\to TA$ can be identified with the
distinguished triangle induced from the short exact sequence $0\to
\tilde{H}A\to \tilde{H}B\to \tilde{H}C\to 0$, which implies the
result.

General case: There is an octahedral diagram
\begin{equation*}
\xymatrix{ A \ar[r] \ar@{=}[d] & \tau_{\leq i}B \ar[r] \ar[d] &
\tau_{\leq i}C
\ar[r] \ar[d] & TA \ar@{=}[d] \\
A \ar[r] & B \ar[r] \ar[d] & C \ar[r] \ar[d] & TA \\
& \tau_{>i}B \ar@{=}[r] \ar[d] & \tau_{>i}B \ar[d] & \\
& T\tau_{\leq i}B \ar[r] & T\tau_{\leq i}C & }
\end{equation*}
and a corresponding exact diagram of cohomology objects. It follows
that the front face and back face in the following diagram are
commutative.
\begin{equation*}
\xymatrix@C-0.5cm{ & {\scriptstyle [\tilde HB]}
\ar[rr]^{f'_2(\Delta)}
\ar'[d][dd] & & {\scriptstyle [\tilde HA]\otimes [\tilde HC]} \ar[dd] \\
{\scriptstyle [B]} \ar[ru]_{\eta_B} \ar[rr]^(0.65){f_2(\Delta)}
\ar[dd] & & {\scriptstyle [A]\otimes [C]}
\ar[ru]_{\eta_A\otimes\eta_C} \ar[dd] & \\
& {\scriptstyle [\tilde H\tau_{\leq i}B]\otimes [\tilde
H\tau_{>i}B]} \ar'[r][rr] & & {\scriptstyle [\tilde HA]\otimes
[\tilde H\tau_{\leq i}C]\otimes [\tilde H\tau_{>i}B]} \\
{\scriptstyle [\tau_{\leq i}B]\otimes [\tau_{>i}B]}
\ar[ru]_{\eta_{\tau_{\leq i}B}\otimes\eta_{\tau_{>i}B}} \ar[rr] & &
{\scriptstyle [A]\otimes [\tau_{\leq i}C]\otimes [\tau_{>i}B]}
\ar[ru]_{\eta_A\otimes\eta_{\tau_{\leq i}C}\otimes\eta_{\tau_{>i}B}}
& }
\end{equation*}
The left and right hand faces are commutative by case (i) considered
above (note that $\tau_{>i}B\cong\tau_{>i}C$ canonically), and the
bottom face is commutative by case (ii). It follows that the top
face is commutative which is the statement of the lemma.
\end{proof}

For any object $A$ in $\calT$ we let $\mu_A$ denote the unit
structure on $[A]\otimes [TA]$ which is induced by the determinant
functor $f$ (compare Lemma \ref{lemma_mu}). Using the canonical
isomorphism $TJ\tilde HA\cong JT\tilde HA$, we will consider
$\mu_{J\tilde HA}$ also as a unit structure on $[J\tilde HA]\otimes
[JT\tilde HA]=[\tilde HA]\otimes [T\tilde HA]$.

\begin{lemma}
\label{lemma_cohomology_2} Let $A$ be an object in $\calT$. Then the
composite map $[A]\otimes [TA]\xrightarrow{\eta_A\otimes\eta_{TA}}
[J\tilde HA]\otimes [J\tilde HTA]\xrightarrow{\cong} [J\tilde
HA]\otimes [JT\tilde HA]$ (where we use the canonical isomorphism
$J\tilde HTA\cong JT\tilde HA$) is an isomorphism of units
$([A]\otimes [TA],\mu_A)\to ([\tilde HA]\otimes [T\tilde
HA],\mu_{J\tilde HA})$.
\end{lemma}

\begin{proof}
If $A\in\calT^{\geq i}\cap\calT^{\leq i}$ then the statement is easy
to see. In general we use induction on the number of non-zero
cohomology objects of $A$, where the inductive step uses Lemmas
\ref{lemma_cohomology_1} and \ref{lemma_mu}(ii) applied to the
distinguished triangles $\tau_{\leq i}A\xrightarrow{a}
A\xrightarrow{b} \tau_{>i}A\xrightarrow{c} T\tau_{\leq i}A$ and
$T\tau_{\leq i}A\xrightarrow{Ta} TA\xrightarrow{Tb}
T\tau_{>i}A\xrightarrow{-Tc} T^2\tau_{\leq i}A$.
\end{proof}

\begin{lemma}
\label{lemma_cohomology_3} Suppose that the diagram of the form
(\ref{diagram_cohomology}) associated to the rotated triangle
$\Delta': B\xrightarrow{b} C\xrightarrow{c}TA\xrightarrow{-Ta}TB$ is
commutative. Then the diagram (\ref{diagram_cohomology}) associated
to the original triangle $\Delta: A\xrightarrow{a} B\xrightarrow{b}
C\xrightarrow{c} TA$ is commutative.
\end{lemma}

\begin{proof}
Consider the following diagram.
\begin{equation*}
\xymatrix{ {\scriptstyle [B]} \ar[r]^-{\mu_{A}} \ar[d]^{\eta_B} &
{\begin{array}{c} \scriptstyle [A]\otimes[TA] \\
\scriptstyle \otimes[B] \end{array}} \ar[rr]^{f_2(\Delta')^{-1}}
\ar[d]^{\eta_A\otimes\eta_{TA}\otimes\eta_B}
& & {\scriptstyle [A]\otimes[C]} \ar[d]^{\eta_A\otimes\eta_C} \\
{\scriptstyle [\tilde HB]} \ar[r]^-{\mu_{J\tilde HA}} \ar@{=}[dd] &
{\begin{array}{c} \scriptstyle [\tilde HA]\otimes[T\tilde HA] \\
\scriptstyle \otimes[\tilde HB] \end{array}} \ar[r]^-{\alpha_1}
\ar[d]^-{\beta_1} &
{\begin{array}{c} \scriptstyle [\tilde HA]\otimes[T\ker \tilde Ha]\otimes[T\ker \tilde Hb] \\
\scriptstyle \otimes [\ker \tilde Hb]\otimes[\ker \tilde Hc]
\end{array}} \ar[r]^-{\alpha_2} \ar[d]^-{\beta_2} & {\scriptstyle
[\tilde HA]\otimes[\tilde HC]} \ar@{=}[dd] \\
& {\begin{array}{c} \scriptstyle [\ker \tilde Ha]\otimes[T\ker \tilde Ha] \\
\scriptstyle \otimes[\ker \tilde Hb]\otimes[T\ker \tilde Hb] \\
\scriptstyle \otimes[\tilde HB]
\end{array}} \ar[d]^-{\mu_{J\ker \tilde Hb}^{-1}} \ar[r]^{\alpha_3} &
{\begin{array}{c} \scriptstyle [\ker \tilde Ha]\otimes[T\ker \tilde Ha] \\
\scriptstyle \otimes[\ker \tilde Hb]\otimes[T\ker \tilde Hb] \\
\scriptstyle \otimes[\ker \tilde Hb]\otimes[\ker \tilde Hc]
\end{array}} \ar[d]^-{\mu_{J\ker \tilde Hb, \mathrm{top}}^{-1}}
& \\
{\scriptstyle [\tilde HB]} \ar[r]^-{\mu_{J\ker \tilde Ha}} &
{\begin{array}{c} \scriptstyle [\ker \tilde Ha]\otimes[T\ker \tilde Ha] \\
\scriptstyle \otimes[\tilde HB]
\end{array}} \ar[r]^-{\alpha_4} &
{\begin{array}{c} \scriptstyle [\ker \tilde Ha]\otimes[T\ker \tilde Ha] \\
\scriptstyle \otimes[\ker \tilde Hb]\otimes[\ker \tilde Hc]
\end{array}} \ar[r]^-{\alpha_5} & {\scriptstyle [\tilde HA]\otimes[\tilde HC]} }
\end{equation*}

We will say that a map in the diagram is induced by a short exact
sequence $\Sigma$ in $\Gr^\rmb(\calE)$ if it is obtained by applying
either $f_2(J(\Sigma))$ or $f_2(J(\Sigma))^{-1}$.

The map $\alpha_1$ is induced by the short exact sequences $\ker
\tilde Hb\to \tilde HB\xrightarrow{\tilde Hb}\ker \tilde Hc$ and
$T\ker \tilde Ha\to T\tilde HA\xrightarrow{-T\tilde Ha}T\ker \tilde
Hb$, and the map $\alpha_2$ is induced by $\mu_{J\ker \tilde
Hb}^{-1}$ and the short exact sequence $\ker \tilde Hc\to \tilde
HC\xrightarrow{\tilde Hc} T\ker \tilde Ha$. The maps $\alpha_3$ and
$\alpha_4$ are induced by the short exact sequence $\ker \tilde
Hb\to \tilde HB\xrightarrow{\tilde Hb}\ker \tilde Hc$, and the map
$\alpha_5$ is induced by the short exact sequences $\ker \tilde
Ha\to \tilde HA\xrightarrow{\tilde Ha}\ker \tilde Hb$ and $\ker
\tilde Hc\to \tilde HC\xrightarrow{\tilde Hc} T\ker \tilde Ha$.

The map $\beta_1$ is induced by the short exact sequences $\ker
\tilde Ha\to \tilde HA\xrightarrow{\tilde Ha}\ker \tilde Hb$ and
$T\ker \tilde Ha\to T\tilde HA\xrightarrow{T\tilde Ha} T\ker \tilde
Hb$. Finally $\beta_2$ is induced by $\ker \tilde Ha\to \tilde
HA\xrightarrow{\tilde Ha}\ker \tilde Hb$ and $\varepsilon([\ker
\tilde Hb])$.

The top left hand square is commutative by Lemma
\ref{lemma_cohomology_2}, and the top right hand rectangle by our
assumption that the statement of the theorem holds for $\Delta'$.
The bottom left hand rectangle is commutative because $\mu$ is
compatible with distinguished triangles, cf.\ Lemma
\ref{lemma_mu}(ii). The two squares in the middle are obviously
commutative (for the upper one note that the minus sign in
$\alpha_1$ and the $\varepsilon([\ker \tilde Hb])$ in $\beta_2$
cancel). Finally the bottom right hand rectangle is commutative
because we are applying $\mu_{J\ker \tilde Hb}$ to two different
copies of $[\ker \tilde Hb]$ which cancels the $\varepsilon([\ker
\tilde Hb])$ in $\beta_2$.

The composite of the two maps in the top row is $f_2(\Delta)$ (by
Lemma \ref{lemma_rotation}) and the composite of the three maps in
the bottom row is the map $f'_2(\Delta): [\tilde HB]\to[\tilde
HA]\otimes[\tilde HC]$ induced by the long exact cohomology
sequence. Hence the commutativity of the outer square is the
required result.
\end{proof}

We can now show that diagram (\ref{diagram_cohomology}) is
commutative for every distinguished triangle $\Delta:
A\xrightarrow{a} B\xrightarrow{b} C\xrightarrow{c} TA$. We do this
by induction on the length of the long exact cohomology sequence of
$\Delta$, where the length is defined as follows: if we write the
cohomology sequence in the form $\dots\to E^{-1}\to E^{0}\to
E^{1}\to E^{2}\to\dots$ where $E^{-1}=H^{-1}(C)$, $E^{0}=H^0(A)$,
$E^{1}=H^0(B)$, $E^{2}=H^0(C)$ etc., then we define the length to be
$0$ if all $E^i$ are zero objects, and to be $j-i$ if $i$ (resp.\
$j$) is the smallest (resp.\ largest) index for which $E^i$ (resp.\
$E^j$) is non-zero.

If the cohomology sequence of $\Delta$ has length $0$ then the
commutativity of (\ref{diagram_cohomology}) is a special case of
Lemma \ref{lemma_cohomology_1}.

Now assume that the cohomology sequence of $\Delta$  contains a
non-zero object. Lemma \ref{lemma_cohomology_3} shows that we can
assume that the lowest non-zero term in the cohomology sequence is
$H^i(A)$ for some $i$. We define an object $D$ of $\calT$ by the
octahedral diagram
\begin{equation}
\label{diagram_def_D}
\begin{split}
\xymatrix{ \tau_{\leq i}A \ar[r] \ar@{=}[d] & A \ar[r] \ar[d]^{a} &
\tau_{> i}A
\ar[r] \ar[d]^{d} & T\tau_{\leq i}A \ar@{=}[d] \\
\tau_{\leq i}A \ar[r] & B \ar[r] \ar[d] & D \ar[r] \ar[d] &
T\tau_{\leq i}A \\
& C \ar@{=}[r] \ar[d] & C \ar[d] & \\
& TA \ar[r] & T\tau_{> i}A & }
\end{split}
\end{equation}
where the first vertical triangle is $\Delta$ and the first
horizontal triangle is the canonical one. We denote the second
vertical triangle by $\Delta'$ and note that the cohomology sequence
of $\Delta'$ is shorter than the cohomology sequence of $\Delta$.
From (\ref{diagram_def_D}) we obtain the following commutative exact
diagram in $\Gr^\rmb(\calE)$.
\begin{equation*}
\xymatrix{
& & 0 \ar[d] & 0 \ar[d] & \\
& & \ker \tilde Ha \ar@{=}[r] \ar[d] & \ker \tilde Hd \ar[d] & \\
0 \ar[r] & T^{-i}H^iA \ar[r] \ar@{=}[d] & \tilde HA \ar[r] \ar[d] &
\tilde H\tau_{>i}A \ar[r] \ar[d] & 0 \\
0 \ar[r] & T^{-i}H^iA \ar[r] & \tilde HB \ar[r] \ar[d] & \tilde HD
\ar[r] \ar[d] & 0 \\
& & \tilde HC \ar@{=}[r] \ar[d] & \tilde HC \ar[d] & \\
& & T\ker \tilde Ha \ar@{=}[r] \ar[d] & T\ker \tilde Hd \ar[d] & \\
& & 0 & 0 & }
\end{equation*}

Finally consider the following diagram in $\calP$.
\begin{equation*}
\xymatrix@C-0.1cm{ & {\scriptstyle [\tilde HB]}
\ar[rr]^{f'_2(\Delta)} \ar'[d][dd] & &
{\scriptstyle [\tilde HA]\otimes [\tilde HC]} \ar[dd] \\
{\scriptstyle [B]} \ar[ur]_{\eta_B} \ar[rr]^(0.65){f_2(\Delta)}
\ar[dd] & & {\scriptstyle [A]\otimes [C]}
\ar[ur]_{\eta_A\otimes\eta_C}
\ar[dd] & \\
& {\scriptstyle [T^{-i}H^iA]\otimes [\tilde HD]}
\ar'[r]^-{\id\otimes f'_2(\Delta')}[rr] & & {\scriptstyle
[T^{-i}H^iA]\otimes [\tilde
H\tau_{>i}A]\otimes [\tilde HC]} \\
{\scriptstyle [\tau_{\leq i}A]\otimes [D]}
\ar[ru]_{\id\otimes\eta_D} \ar[rr]^{\id\otimes f_2(\Delta')} &&
{\scriptstyle [\tau_{\leq i}A]\otimes [\tau_{>i}A]\otimes [C]}
\ar[ru]_{\id\otimes\eta_{\tau_{>i}A}\otimes\eta_C} & }
\end{equation*}
The commutativity of the bottom face follows from the induction
hypothesis applied to $\Delta': \tau_{> i}A\to D\to C\to T\tau_{>i}
A$, and the commutativity of the back face follows easily from the
exact commutative diagram above. The front is commutative by the
associativity axiom applied to diagram (\ref{diagram_def_D}). The
commutativity of the left and right hand faces are special cases of
Lemma \ref{lemma_cohomology_1}. The resulting commutativity of the
top face gives the commutativity of diagram
(\ref{diagram_cohomology}) and therefore completes the proof of
Proposition \ref{proposition_cohomology}.
\end{proof}


\end{document}